\documentclass[leqno,draft]{article}

\def\dim{\mathop{\rm dim}}
\def\diam{\mathop{\rm diam}}
\def\dist{\mathop{\rm dist}}

\newtheorem{theorem}{Theorem}
\newtheorem{lemma}[theorem]{Lemma}
\newtheorem{proposition}[theorem]{Proposition}
\newtheorem{definition}[theorem]{Definition}
\newtheorem{corollary}[theorem]{Corollary}

\newcommand{\begintheorem}{\addtocounter{equation}{1}\begin{theorem}}
\newcommand{\beginlemma}{\addtocounter{equation}{1}\begin{lemma}}
\newcommand{\beginproposition}{\addtocounter{equation}{1}\begin{proposition}}
\newcommand{\begindefinition}{\addtocounter{equation}{1}\begin{definition}}
\newcommand{\begincorollary}{\addtocounter{equation}{1}\begin{corollary}}

\begin{document}

\title{Some elementary aspects of \\ Hausdorff measure and dimension}

\author{Stephen Semmes \\
        Rice University}

\date{}

\maketitle

\begin{abstract}
Basic properties of Hausdorff content, measure, and dimension of
subsets of metric spaces are discussed, especially in connection with
Lipschitz mappings and topological dimension.
\end{abstract}

\tableofcontents

\section{Metric spaces}
\label{metric spaces}
\setcounter{equation}{0}

        Let $(M, d(x, y))$ be a metric space.  Thus $M$ is a set, and
$d(x, y)$ is a nonnegative real-valued function defined for $x, y \in
M$ such that $d(x, y) = 0$ if and only if $x = y$,
\begin{equation}
        d(y, x) = d(x, y)
\end{equation}
for every $x, y \in M$, and
\begin{equation}
        d(x, z) \le d(x, y) + d(y, z)
\end{equation}
for every $x, y, z \in M$.

        If $p \in M$ and $r > 0$, then the open ball in $M$ with
center $p$ and radius $r$ is defined by
\begin{equation}
        B(p, r) = \{q \in M : d(p, q) < r\}.
\end{equation}
A set $A \subseteq M$ is said to be \emph{bounded} if it is contained
in a ball, in which case its \emph{diameter} is defined by
\begin{equation}
        \diam A = \sup \{d(x, y) : x, y \in A\}.
\end{equation}
This can be interpreted as being $0$ when $A = \emptyset$, and it is
sometimes convenient to make the convention that the diameter of an
unbounded set is $+ \infty$.  The \emph{closure} $\overline{A}$ of a
set $A \subseteq M$ is defined to be the union of $A$ and the set of
limit points of $A$ in $M$.  It is easy to see that the closure of a
bounded set $A$ is also bounded, and has the same diameter as $A$.

        A metric space $M$ is said to be \emph{separable} if there is
a dense set $E \subseteq M$ with only finitely or countably many
elements.  For example, the real line ${\bf R}$ with the standard
metric is separable, because the set ${\bf Q}$ of rational numbers is
countable and dense in ${\bf R}$.  Similarly, ${\bf R}^n$ is separable
for each positive integer $n$ with respect to the standard Euclidean
metric, because ${\bf Q}^n$ is a countable dense set in ${\bf R}^n$.
Note that a metric space $M$ is separable if and only if for each $r >
0$ there is a finite or countable set $E_r \subseteq M$ which is
$r$-dense in $M$ in the sense that
\begin{equation}
\label{M = bigcup_{y in E_r} B(y, r)}
        M = \bigcup_{y \in E_r} B(y, r).
\end{equation}
Equivalently, (\ref{M = bigcup_{y in E_r} B(y, r)}) says that for
every $x \in M$ there is a $y \in E_r$ which satisfies $d(x, y) < r$.
For if $M$ is separable, and $E \subseteq M$ is a dense set with only
finitely or countably many elements, then one can take $E_r = E$ for
every $r > 0$.  Conversely, if for every $r > 0$ there is an $r$-dense
set $E_r$ in $M$ with only finitely or countably many elements, then
\begin{equation}
        E = \bigcup_{n = 1}^\infty E_{1/n}
\end{equation}
is a dense set in $M$ with only finitely or countably many elements,
and hence $M$ is separable.

        A set $A \subseteq M$ is \emph{totally bounded} if for every
$r > 0$ there finitely many elements $x_1, \ldots, x_n$ of $A$ such that
\begin{equation}
        A \subseteq \bigcup_{i = 1}^n B(x_i, r).
\end{equation}
Totally bounded sets are bounded, and bounded subsets of ${\bf R}^n$
are totally bounded.  Totally bounded metric spaces are automatically
separable, by the remarks of the preceding paragraph.  Compact sets
are totally bounded, and it is well known that totally bounded closed
subsets of complete metric spaces are compact.

\section{Hausdorff content, 1}
\label{hausdorff content, 1}
\setcounter{equation}{0}

        Let $(M, d(x, y))$ be a metric space, and let $\alpha$ be a
positive real number.  If $A \subseteq M$, then put
\begin{equation}
 \widetilde{H}^\alpha_{con}(A) = \inf \bigg\{\sum_{i = 1}^n (\diam E_i)^\alpha
  : E_1, \ldots, E_n \subseteq M, \ A \subseteq \bigcup_{i = 1}^n E_i\bigg\}.
\end{equation}
More precisely, we take the infimum of the sum
\begin{equation}
\label{sum_{i = 1}^n (diam E_i)^alpha}
        \sum_{i = 1}^n (\diam E_i)^\alpha
\end{equation}
over all coverings of $A$ by finitely many subsets $E_1, \ldots, E_n$
of $A$ in $M$, where $n$ depends on the covering.  If any of the
$E_i$'s is unbounded, so that $\diam E_i = +\infty$, then (\ref{sum_{i
= 1}^n (diam E_i)^alpha}) is $+\infty$ too.  Otherwise, we can
restrict our attention to coverings of $A$ by finitely many bounded
subsets of $M$ when $A$ is bounded, and set
\begin{equation}
        \widetilde{H}^\alpha_{con}(A) = +\infty
\end{equation}
when $A$ is unbounded.

        If $\alpha = 0$, then let us make the following conventions.
First, $(\diam E)^0 = 0$ when $E = \emptyset$.  Second, $(\diam E)^0 =
1$ when $E \ne \emptyset$ and $E$ is bounded.  Third, $(\diam E)^0 =
+\infty$ when $E$ is unbounded.  In this way, we can extend the
definition of $\widetilde{H}^\alpha_{con}(A)$ to $\alpha = 0$.

        For every $A \subseteq M$,
\begin{equation}
\label{widetilde{H}^alpha_{con}(A) le (diam A)^alpha}
        \widetilde{H}^\alpha_{con}(A) \le (\diam A)^\alpha
\end{equation}
This basically corresponds to the $n = 1$ case of the definition of
$\widetilde{H}^\alpha_{con}(A)$, which is to say that $A$ is covered
by only one set.  By allowing coverings by several sets,
$\widetilde{H}^\alpha_{con}(A)$ may be reduced significantly.  For
instance, if $A$ is the union of two or more subsets of $M$ with small
diameter which are relatively far from each other and $\alpha > 0$,
then $\widetilde{H}^\alpha_{con}(A)$ is much smaller than $(\diam
A)^\alpha$.  Coverings by several sets is not important in the
definition of $\widetilde{H}^\alpha_{con}(A)$ when $\alpha = 0$, in
which event we have equality in (\ref{widetilde{H}^alpha_{con}(A) le
(diam A)^alpha}).

        Observe that $\widetilde{H}^\alpha_{con}(A)$ is monotone in
$A$, in the sense that
\begin{equation}
\label{widetilde{H}^alpha_{con}(A) le widetilde{H}^alpha_{con}(B)}
        \widetilde{H}^\alpha_{con}(A) \le \widetilde{H}^\alpha_{con}(B)
\end{equation}
when $A \subseteq B \subseteq M$.  This is because every covering of
$B$ is also a covering of $A$.  In particular,
\begin{equation}
        \widetilde{H}^\alpha_{con}(B) = 0
\end{equation}
implies that
\begin{equation}
        \widetilde{H}^\alpha_{con}(A) = 0
\end{equation}
when $A \subseteq B$.  Of course,
$\widetilde{H}^\alpha_{con}(\emptyset) = 0$.

\section{Closed sets}
\label{closed sets}
\setcounter{equation}{0}

        If we restrict our attention to coverings of $A \subseteq M$
by finitely many closed subsets $E_1, \ldots, E_n$ of $M$, then
$\widetilde{H}^\alpha_{con}(A)$ would still be the same as before.
For if $E_1, \ldots, E_n$ are arbitrary subsets of $M$ such that
\begin{equation}
        A \subseteq \bigcup_{i = 1}^n E_i,
\end{equation}
then their closures $\overline{E_1}, \ldots, \overline{E_n}$
are closed subsets of $M$ which satisfy
\begin{equation}
        A \subseteq \bigcup_{i = 1}^n \overline{E_i}
\end{equation}
and
\begin{equation}
        \sum_{i = 1}^n (\diam \overline{E_i})^\alpha
          = \sum_{i = 1}^n (\diam E_i)^\alpha,
\end{equation}
since $\diam \overline{E_i} = \diam E_i$ for each $i$.  As a
consequence,
\begin{equation}
\label{widetilde{H}^alpha_{con}(overline{A}) = widetilde{H}^alpha_{con}(A)}
 \widetilde{H}^\alpha_{con}(\overline{A}) = \widetilde{H}^\alpha_{con}(A)
\end{equation}
for any set $A \subseteq M$.  Indeed, any covering of $A$ by finitely
many closed sets also covers $\overline{A}$, since the union of
finitely many closed sets is closed.

\section{Open sets}
\label{open sets}
\setcounter{equation}{0}

        If $E \subseteq M$ and $r > 0$, then put
\begin{equation}
        E(r) = \bigcup_{x \in E} B(x, r).
\end{equation}
This is the same as the set of $y \in M$ for which there is an $x \in
E$ such that $d(x, y) < r$.  In particular, $E(r) = M$ exactly when
$E$ is $r$-dense in $M$.

        By construction, $E(r)$ is an open set in $M$, and
\begin{equation}
        E \subseteq E(r).
\end{equation}
If $E$ is bounded, then $E(r)$ is bounded, and
\begin{equation}
        \diam E(r) \le \diam E + 2 \, r.
\end{equation}
This is easy to see, directly from the definitions.

        It follows that $\widetilde{H}^\alpha_{con}(A)$ would also be
the same if we restricted ourselves to coverings of $A$ by finitely many
open subsets of $M$.  For if $E_1, \ldots, E_n$ is any collection of
finitely many subsets of $M$ that covers $A$, then $E_1(r), \ldots, E_n(r)$
are finitely many open subsets of $M$ covering $A$ for each $r > 0$, and
\begin{equation}
        \lim_{r \to 0} \sum_{i = 1}^n (\diam E_i(r))^\alpha
          = \sum_{i = 1}^n (\diam E_i)^\alpha.
\end{equation}

\section{Intervals}
\label{intervals}
\setcounter{equation}{0}

        If $a$, $b$ are real numbers with $a < b$, then the open and
closed intervals in the real line with endpoints $a$, $b$ are defined
by
\begin{equation}
        (a, b) = \{x \in {\bf R} : a < x < b\}
\end{equation}
and
\begin{equation}
        [a, b] = \{x \in {\bf R} : a \le x \le b\},
\end{equation}
respectively.  The latter also makes sense when $a = b$, and reduces
to a single point.  Note that
\begin{equation}
        \diam \, (a, b) = \diam \, [a, b] = b - a.
\end{equation}

        If $E$ is a bounded nonempty set in the real line, then
\begin{equation}
        I = [\inf E, \sup E]
\end{equation}
is a closed interval that contains $E$ and has the same diameter.  If
$E$ is unbounded, then one can consider ${\bf R}$ as an unbounded
interval that contains $E$.  Similarly, any set $E \subseteq {\bf R}$
is contained in an open interval whose diameter is arbitrarily close
to the diameter of $E$.

        As in the previous sections, it follows that
$\widetilde{H}^\alpha_{con}(A)$ would be the same when $M = {\bf R}$
with the standard metric if we restricted our attention to coverings
of $A$ by finitely many open or closed intervals.

        It is well known that
\begin{equation}
        \widetilde{H}^1_{con}((a, b)) = \widetilde{H}^1_{con}([a, b]) = b - a.
\end{equation}
By (\ref{widetilde{H}^alpha_{con}(A) le (diam A)^alpha}) and
(\ref{widetilde{H}^alpha_{con}(A) le widetilde{H}^alpha_{con}(B)}),
\begin{equation}
 \widetilde{H}^1_{con}((a, b)) \le \widetilde{H}^1_{con}([a, b]) \le b - a,
\end{equation}
and so one only has to show that
\begin{equation}
        b - a \le \widetilde{H}^1_{con}((a, b)).
\end{equation}
To see this, it suffices to check that
\begin{equation}
        b - a \le \sum_{j = 1}^n \diam I_j
\end{equation}
when $I_1, \ldots, I_n$ are finitely many intervals covering $(a, b)$.

\section{Subadditivity}
\label{subadditivity}
\setcounter{equation}{0}

        In any metric space $M$,
\begin{equation}
        \widetilde{H}^\alpha_{con}(A \cup B)
 \le \widetilde{H}^\alpha_{con}(A) + \widetilde{H}^\alpha_{con}(B)
\end{equation}
for every $A, B \subseteq M$, basically because coverings of $A$ and
$B$ can be combined to give coverings of $A \cup B$.  In particular,
\begin{equation}
        \widetilde{H}^\alpha_{con}(A \cup B) = 0
\end{equation}
when $\widetilde{H}^\alpha_{con}(A) = \widetilde{H}^\alpha_{con}(B) = 0$.
Similarly,
\begin{equation}
 \widetilde{H}^\alpha_{con}(A \cup B) = \widetilde{H}^\alpha_{con}(A)
\end{equation}
for every $A, B \subseteq M$ with $\widetilde{H}^\alpha_{con}(B) = 0$,
because
\begin{equation}
 \widetilde{H}^\alpha_{con}(A) \le \widetilde{H}^\alpha_{con}(A \cup B)
 \le \widetilde{H}^\alpha_{con}(A) + \widetilde{H}^\alpha_{con}(B)
   = \widetilde{H}^\alpha_{con}(A).
\end{equation}
Note that $\widetilde{H}^\alpha_{con}(B) = 0$ when $B$ has only
finitely many elements and $\alpha > 0$, while
$\widetilde{H}^0_{con}(B) \ge 1$ when $B \ne \emptyset$.

\section{Hausdorff content, 2}
\label{hausdorff content, 2}
\setcounter{equation}{0}

        Let $(M, d(x, y))$ be a metric space, and let $\alpha$ be a
nonnegative real number.  If $A \subseteq M$, then put
\begin{equation}
 H^\alpha_{con}(A) = \inf \bigg\{\sum_i (\diam E_i)^\alpha : E_i \subseteq M, 
                       \ A \subseteq \bigcup_i E_i\bigg\},
\end{equation}
where more precisely the infimum is taken over coverings of $A$ by
finitely or countably many subsets $E_i$ of $M$.  In the case of
countable coverings of $A$,
\begin{equation}
        \sum_i (\diam E_i)^\alpha
\end{equation}
is interpreted as the supremum of the subsums of $(\diam E_i)^\alpha$
over finitely many indices $i$, which may be infinite.  We also use
the same conventions for $\alpha = 0$ as before.

        Note that
\begin{equation}
        H^\alpha_{con}(A) \le \widetilde{H}^\alpha_{con}(A),
\end{equation}
since the finite coverings used in the definition of
$\widetilde{H}^\alpha_{con}(A)$ are also included in the definition of
$H^\alpha_{con}(A)$.  In particular,
\begin{equation}
        H^\alpha_{con}(A) \le (\diam A)^\alpha,
\end{equation}
as in (\ref{widetilde{H}^alpha_{con}(A) le (diam A)^alpha}).
We also have the monotonicity property
\begin{equation}
        H^\alpha_{con}(A) \le H^\alpha_{con}(B)
\end{equation}
when $A \subseteq B \subseteq M$, as in the previous case.

        If we restrict our attention to coverings of $A$ by finitely
or countably many closed subsets $E_i$ of $M$, then
$H^\alpha_{con}(A)$ would still be the same, as in Section \ref{closed
sets}.  However, this does not imply that
\begin{equation}
        H^\alpha_{con}(\overline{A}) = H^\alpha_{con}(A),
\end{equation}
since the union of infinitely many closed subsets of $M$ may not be closed.
One can also check that $H^\alpha_{con}(A)$ would be the same if we restricted
our attention to coverings of $A$ by finitely or countably many open subsets
$E_i$ of $M$, as in Section \ref{open sets}.  For if $\{E_i\}_i$
is any collection of finitely or countably many subsets of $M$ that
covers $A$, then one can choose positive real numbers $r_i$ such that
\begin{equation}
        \sum_i (\diam E_i(r_i))^\alpha
\end{equation}
is arbitrarily close to $\sum_i (\diam E_i)^\alpha$.  If $M$ is the
real line with the standard metric, then we can use coverings of $A
\subseteq {\bf R}$ by finitely or countably many open or closed
intervals and get the same result for $H^\alpha_{con}(A)$, as in
Section \ref{intervals}.

        If $K$ is a compact set in any metric space $M$, then
\begin{equation}
        H^\alpha_{con}(K) = \widetilde{H}^\alpha_{con}(K).
\end{equation}
Indeed, we can restrict our attention to coverings of $K$ by open
subsets of $M$, as in the previous paragraph, and then reduce to
finite coverings by compactness.

\section{Countable subadditivity}
\label{countable subadditivity}
\setcounter{equation}{0}

        As before,
\begin{equation}
        H^\alpha_{con}(A \cup B) \le H^\alpha_{con}(A) + H^\alpha_{con}(B)
\end{equation}
for every $A, B \subseteq M$.  Moreover,
\begin{equation}
\label{H^alpha_{con}(cup_j A_j) le sum_j H^alpha_{con}(A_j)}
        H^\alpha_{con}\Big(\bigcup_{j = 1}^\infty A_j\Big)
         \le \sum_{j = 1}^\infty H^\alpha_{con}(A_j)
\end{equation}
for any sequence $A_1, A_2, \ldots$ of subsets of $M$.  Of course, the
sum on the right may be $+\infty$, either because $H^\alpha_{con}(A_j)
= +\infty$ for some $j$, or $H^\alpha_{con}(A_j) < + \infty$ for each
$j$ but the sum diverges, in which event the inequality is trivial.

        To prove (\ref{H^alpha_{con}(cup_j A_j) le sum_j
H^alpha_{con}(A_j)}), we may as well suppose that $H^\alpha_{con}(A_j)
< + \infty$ for each $j$, and that their sum converges.  Let $\epsilon
> 0$ be given, and for each $j$, let $\{E_{j, l}\}_l$ be a collection
of finitely or countably many subsets of $M$ such that
\begin{equation}
        A_j \subseteq \bigcup_l E_{j, l}
\end{equation}
and
\begin{equation}
        \sum_l (\diam E_{j, l})^\alpha
          < H^\alpha_{con}(A_j) + 2^{-j} \, \epsilon.
\end{equation}
By combining these families, we get a collection $\{E_{j, l}\}_{j, l}$
of finitely or countably many subsets of $M$ such that
\begin{equation}
        \bigcup_{j = 1}^\infty A_j \subseteq \bigcup_{j, l} E_{j, l}
\end{equation}
and
\begin{equation}
        \sum_{j, l} (\diam E_{j, l})^\alpha
         < \sum_{j = 1}^\infty (H^\alpha_{con}(A_j) + 2^{-j} \, \epsilon)
         = \sum_{j = 1}^\infty H^\alpha_{con}(A_j) + \epsilon.
\end{equation}
Thus
\begin{equation}
        H^\alpha_{con}\Big(\bigcup_{j = 1}^\infty A_j\Big)
         < \sum_{j = 1}^\infty H^\alpha_{con}(A_j) + \epsilon,
\end{equation}
which implies (\ref{H^alpha_{con}(cup_j A_j) le sum_j
H^alpha_{con}(A_j)}), since $\epsilon > 0$ is arbitrary.

        In particular,
\begin{equation}
        H^\alpha_{con}\Big(\bigcup_{j = 1}^\infty A_j\Big) = 0
\end{equation}
when $H^\alpha_{con}(A_j) = 0$ for each $j$.  It follows that
\begin{equation}
        H^\alpha_{con}(A) = 0
\end{equation}
for every countable set $A$ when $\alpha > 0$.

        By contrast, if $M$ is unbounded, then there are unbounded
countable subsets $A$ of $M$, for which
\begin{equation}
        \widetilde{H}^\alpha_{con}(A) = +\infty.
\end{equation}
If $M$ is the real line with the standard metric, and $A$ is the set
of rational numbers in an interval $[a, b]$, $a < b$, then $A$ is a
bounded countable set, while
\begin{equation}
        \widetilde{H}^1_{con}(A) = \widetilde{H}^1_{con}(\overline{A})
         = \widetilde{H}^1_{con}([a, b]) = b - a > 0.
\end{equation}
Observe that
\begin{equation}
        \widetilde{H}^1_{con}(\overline{A}) = H^1_{con}(\overline{A}),
\end{equation}
since $\overline{A} = [a, b]$ is compact, and so we get examples where
\begin{equation}
        H^1_{con}(A) < H^1_{con}(\overline{A}),
\end{equation}
as in the previous section.

\section{Borel measures}
\label{borel measures}
\setcounter{equation}{0}

        Let $\mu$ be a nonnegative Borel measure on a metric space $M$
such that
\begin{equation}
\label{mu(E) le C (diam E)^alpha}
        \mu(E) \le C \, (\diam E)^\alpha
\end{equation}
for some $\alpha, C > 0$ and every Borel set $E$ in $M$.  For example,
Lebesgue measure on ${\bf R}^n$ has this property with $\alpha = n$,
and with $C = 1$ when $n = 1$.  In this case,
\begin{equation}
        \mu(A) \le C \, H^\alpha_{con}(A)
\end{equation}
for every Borel set $A$ in $M$.  For if $\{E_i\}_i$ is any collection
of finitely or countably many Borel sets in $M$ such that
\begin{equation}
        A \subseteq \bigcup_i E_i,
\end{equation}
then
\begin{equation}
        \mu(A) \le \sum_i \mu(E_i) \le C \, \sum_i (\diam E_i)^\alpha.
\end{equation}
The desired estimate follows, since it suffices to consider coverings
of $A$ by open or closed subsets of $M$ in the definition of
$H^\alpha_{con}(A)$, which are automatically Borel sets.  This is a
basic technique to get lower bounds for $H^\alpha_{con}(A)$, while
upper bounds can be obtained from explicit coverings of $A$.  In
particular, this is a common way to show that $H^\alpha_{con}(A) > 0$.

\section{Localization}
\label{localization}
\setcounter{equation}{0}

        Let $M$ be a metric space, and fix $\alpha \ge 0$, $0 < \delta
\le \infty$.  If $A \subseteq M$, then put
\begin{eqnarray}
 \widetilde{H}^\alpha_\delta(A) & = & 
  \inf \bigg\{\sum_{i = 1}^n (\diam E_i)^\alpha : E_1, \ldots, E_n \subseteq M,
     \ A \subseteq \bigcup_{i = 1}^n E_i,  \\
 & & \qquad\qquad\qquad\quad \hbox{and } \diam E_i < \delta
             \hbox{ for } i = 1, \ldots, n\bigg\} \nonumber
\end{eqnarray}
and
\begin{eqnarray}
 H^\alpha_\delta(A) & = & \inf\bigg\{\sum_i (\diam E_i)^\alpha :
     \ E_i \subseteq M, \ A \subseteq \bigcup_i E_i, \hbox{ and} \\
& & \qquad\qquad\qquad\qquad\quad \diam E_i < \delta
                      \hbox{ for each } i\bigg\}. \nonumber
\end{eqnarray}
In the first case, the infimum is taken over all coverings of $A$ by
finitely many sets with diameter less than $\delta$, while in the
second case the infimum is taken over all coverings of $A$ by finitely
or countably many sets with diameter less than $\delta$.  This is
interpreted as being $+\infty$ when there is no such covering of $A$.
If $A$ is totally bounded, then there are admissible coverings of $A$
for $\widetilde{H}^\alpha_\delta(A)$ for every $\delta > 0$.
Similarly, there are admissible coverings of $A$ for
$H^\alpha_\delta(A)$ for every $\delta > 0$ when $M$ is separable.
However, the corresponding infinite sums $\sum_i (\diam E_i)^\alpha$
may still diverge, so that $H^\alpha_\delta(A) = +\infty$.  Of course,
\begin{equation}
        H^\alpha_\delta(A) \le \widetilde{H}^\alpha_\delta(A),
\end{equation}
since the coverings of $A$ used in the definition of
$\widetilde{H}^\alpha_\delta(A)$ are more restrictive than the
coverings used in the definition of $H^\alpha_\delta(A)$.

        By definition,
\begin{equation}
        \widetilde{H}^\alpha_{con}(A) = \widetilde{H}^\alpha_\infty(A)
\end{equation}
and
\begin{equation}
        H^\alpha_{con}(A) = H^\alpha_\infty(A)
\end{equation}
for every $A \subseteq M$ and $\alpha \ge 0$.  Also,
\begin{equation}
        \widetilde{H}^\alpha_\eta(A) \le \widetilde{H}^\alpha_\delta(A)
\end{equation}
and
\begin{equation}
        H^\alpha_\eta(A) \le H^\alpha_\delta(A)
\end{equation}
when $0 < \delta \le \eta \le \infty$.  This is because the infima in
the definitions of $\widetilde{H}^\alpha_\delta(A)$,
$H^\alpha_\delta(A)$ are taken over more restrictive classes of
coverings of $A$ as $\delta$ decreases, which implies that these
infima are increasing as $\delta$ decreases.  As a partial converse to
this statement, $\widetilde{H}^\alpha_\delta(A)$ or
$H^\alpha_\delta(A)$ is equal to $0$ for every $\delta > 0$ as soon as
it is $0$ for some $\delta > 0$.  For the sets $E_i$ in the
appropriate coverings of $A$ have to have small diameter when the
corresponding sums $\sum_i (\diam E_i)^\alpha$ are small.

        As before,
\begin{equation}
        \widetilde{H}^\alpha_\delta(A) \le \widetilde{H}^\alpha_\delta(B)
\end{equation}
and
\begin{equation}
        H^\alpha_\delta(A) \le H^\alpha_\delta(B)
\end{equation}
when $A \subseteq B$, since every admissible covering of $B$ for
$\widetilde{H}^\alpha_\delta(B)$ or $H^\alpha_\delta(B)$ is also an
admissible covering of $A$.  Moreover,
$\widetilde{H}^\alpha_\delta(A)$ is finitely subadditive and
$H^\alpha_\delta(A)$ is countably subadditive for every $\delta > 0$,
for the same reasons as when $\delta = +\infty$.  We may restrict our
attention to coverings of $A$ by open or closed sets in the
definitions of $\widetilde{H}^\alpha_\delta(A)$ and
$H^\alpha_\delta(A)$, or to coverings by intervals when $M = {\bf R}$,
again for the same reasons as when $\delta = +\infty$.  This implies that
\begin{equation}
 \widetilde{H}^\alpha_\delta(\overline{A}) = \widetilde{H}^\alpha_\delta(A)
\end{equation}
for every $A \subseteq M$, and that
\begin{equation}
        H^\alpha_\delta(K) = \widetilde{H}^\alpha_\delta(K)
\end{equation}
when $K \subseteq M$ is compact, as when $\delta = +\infty$.

        If $A, B \subseteq M$ and
\begin{equation}
        d(x, y) \ge \delta
\end{equation}
for every $x \in A$ and $y \in B$, then
\begin{equation}
        \widetilde{H}^\alpha_\delta(A \cup B) =
 \widetilde{H}^\alpha_\delta(A) + \widetilde{H}^\alpha_\delta(B)
\end{equation}
and
\begin{equation}
        H^\alpha_\delta(A \cup B) = H^\alpha_\delta(A) + H^\alpha_\delta(B).
\end{equation}
For if $\{E_i\}_i$ is a covering of $A \cup B$ by sets of diameter
less than $\delta$, then each $E_i$ will intersect at most one of $A$
and $B$.  This implies that there are disjoint subcollections of
$\{E_i\}_i$ covering $A$ and $B$, which permits one to estimate the sum
of the measures of $A$ and $B$ by the corresponding measure of $A \cup B$.

\section{Hausdorff measures}
\label{hausdorff measures}
\setcounter{equation}{0}

        Let $M$ be a metric space, and fix $\alpha \ge 0$.  If $A
\subseteq M$, then put
\begin{equation}
 \widetilde{H}^\alpha(A) = \sup_{\delta > 0} \widetilde{H}^\alpha_\delta(A)
\end{equation}
and
\begin{equation}
        H^\alpha(A) = \sup_{\delta > 0} H^\alpha_\delta(A).
\end{equation}
The latter is known as the \emph{$\alpha$-dimensional Hausdorff
measure} of $A$.  Note that $\widetilde{H}^\alpha_\delta(A)$ or
$H^\alpha_\delta(A)$ may be $+\infty$ for some $\delta > 0$, or they
may be finite for every $\delta$ but unbounded, so that the supremum
is $+\infty$.  The supremum over $\delta > 0$ can also be considered
as a limit as $\delta \to 0$, since $\widetilde{H}^\alpha_\delta(A)$
and $H^\alpha_\delta(A)$ are monotone increasing as $\delta$
decreases.

        As usual,
\begin{equation}
        H^\alpha(A) \le \widetilde{H}^\alpha(A)
\end{equation}
for every $A \subseteq M$, and
\begin{equation}
        \widetilde{H}^\alpha(A) \le \widetilde{H}^\alpha(B)
\end{equation}
and
\begin{equation}
        H^\alpha(A) \le H^\alpha(B)
\end{equation}
when $A \subseteq B \subseteq M$.  In addition, $\widetilde{H}^\alpha$
is finitely subadditive and $H^\alpha$ is countably subadditive,
because of the corresponding properties of
$\widetilde{H}^\alpha_\delta$ and $H^\alpha_\delta$.  Similarly,
\begin{equation}
        \widetilde{H}^\alpha(\overline{A}) = \widetilde{H}^\alpha(A)
\end{equation}
for every $A \subseteq M$, and
\begin{equation}
        H^\alpha(K) = \widetilde{H}^\alpha(K)
\end{equation}
when $K$ is compact.  Of course,
\begin{equation}
        \widetilde{H}^\alpha_{con}(A) \le \widetilde{H}^\alpha(A)
\end{equation}
and
\begin{equation}
        H^\alpha_{con}(A) \le H^\alpha(A).
\end{equation}
Conversely,
\begin{equation}
 \widetilde{H}^\alpha_{con}(A) = 0 \hbox{ implies } \widetilde{H}^\alpha(A) = 0
\end{equation}
and
\begin{equation}
        H^\alpha_{con}(A) = 0 \hbox{ implies } H^\alpha(A) = 0,
\end{equation}
by the analogous remarks in the previous section.

        If $A, B \subseteq M$ and
\begin{equation}
        d(x, y) \ge \eta
\end{equation}
for some $\eta > 0$ and every $x \in A$, $y \in B$, then
\begin{equation}
        \widetilde{H}^\alpha(A \cup B) = \widetilde{H}^\alpha(A)
                                         + \widetilde{H}^\alpha(B)
\end{equation}
and
\begin{equation}
        H^\alpha(A \cup B) = H^\alpha(A) + H^\alpha(B),
\end{equation}
by the analogous statement in the previous section.  For example, one
can use this to show that $\widetilde{H}^0 = H^0$ is the same as
counting measure.  If $M = {\bf R}$, then $H^1(A)$ is the same as the
Lebesgue outer measure of $A$.  In this case, one can check that
$\widetilde{H}^1_\delta(A)$ and $H^1_\delta(A)$ do not depend on $\delta$.

        If $\alpha < \beta$, then it is easy to see that
\begin{equation}
        \widetilde{H}^\beta_\delta(A)
          \le \delta^{\beta - \alpha} \, \widetilde{H}^\alpha_\delta(A)
\end{equation}
and
\begin{equation}
        H^\beta_\delta(A) \le \delta^{\beta - \alpha} \, H^\alpha_\delta(A).
\end{equation}
This implies that $\widetilde{H}^\beta(A) = 0$ when
$\widetilde{H}^\alpha(A) < +\infty$, and that $H^\beta(A) = 0$ when
$H^\alpha(A) < +\infty$.  Equivalently, $\widetilde{H}^\alpha(A) =
+\infty$ when $\widetilde{H}^\beta(A) > 0$, and $H^\alpha(A) =
+\infty$ when $H^\beta(A) > 0$.

\section{Borel sets}
\label{borel sets}
\setcounter{equation}{0}

        Let $M$ be a metric space, and suppose that $A \subseteq M$ satisfies
\begin{equation}
        H^\alpha(A) < +\infty
\end{equation}
for some $\alpha \ge 0$.  For each positive integer $l$, let $\{E_{j,
l}\}_j$ be a collection of finitely or countably many open subsets of
$M$ such that
\begin{equation}
        \diam E_{j, l} < \frac{1}{l}
\end{equation}
for each $j$,
\begin{equation}
        A \subseteq \bigcup_j E_{j, l},
\end{equation}
and
\begin{equation}
        \sum_j (\diam E_{j, l})^\alpha < H^\alpha_{1/l}(A) + \frac{1}{l}.
\end{equation}
Put
\begin{equation}
        B = \bigcap_{l = 1}^\infty \Big(\bigcup_j E_{j, l}\Big),
\end{equation}
so that $A \subseteq B$ by construction.  Hence $H^\alpha(A) \le
H^\alpha(B)$, while
\begin{equation}
        H^\alpha_{1/l}(B) \le \sum_j (\diam E_{j, l})^\alpha
                             < H^\alpha_{1/l}(A) + \frac{1}{l}
\end{equation}
for each $l$, because $\{E_{j, l}\}_j$ also covers $B$.  Thus
$H^\alpha(B) \le H^\alpha(A)$, and therefore
\begin{equation}
        H^\alpha(A) = H^\alpha(B).
\end{equation}
Note that $\bigcup_j E_{j, l}$ is an open set for each $l$, since the
$E_{j, l}$'s are open sets.  It follows that $B$ is the intersection
of a sequence of open sets, and hence a Borel set.

\section{Borel measures, 2}
\label{Borel measures, 2}
\setcounter{equation}{0}

        Let $M$ be a metric space, let $\mu$ be a Borel measure on
$M$, and fix $\alpha, \delta > 0$.  If
\begin{equation}
\label{mu(E) le C (diam E)^alpha, diam E < delta}
        \mu(E) \le C \, (\diam E)^\alpha
\end{equation}
for some $C > 0$ and every Borel set $E \subseteq M$ with $\diam E <
\delta$, then
\begin{equation}
\label{mu(A) le C H^alpha_delta(A)}
        \mu(A) \le C \, H^\alpha_\delta(A),
\end{equation}
for every Borel set $A \subseteq M$, by the same argument as in
Section \ref{borel measures}.  Actually, it suffices to restrict our
attention to open sets $E$, since we can use coverings by open sets in
the definition of $H^\alpha_\delta$.  However, one can also check that
(\ref{mu(E) le C (diam E)^alpha, diam E < delta}) automatically holds
for arbitrary Borel sets when it holds for open sets and $\mu$ is
outer regular.  It suffices as well to consider only sets $E$ such
that
\begin{equation}
        A \cap E \ne \emptyset.
\end{equation}
This is because sets disjoint from $A$ can be dropped from any
covering of $A$ without affecting the estimates.  Put
\begin{eqnarray}
 X(\alpha, \delta, C) & = & \{x \in M : \mu(E) \le C \, (\diam E)^\alpha
                          \hbox{ for every open set } \\
 & & \qquad\qquad E \subseteq M \hbox{ with } x \in E \hbox{ and }
                           \diam E < \delta\}. \nonumber
\end{eqnarray}
The preceding remarks imply that (\ref{mu(A) le C H^alpha_delta(A)})
holds when $A \subseteq X(\alpha, \delta, C)$.  An advantage of
working with open sets $E$ is that $X(\alpha, \delta, C)$ is
automatically a closed set in $M$, because its complement is an open
set.

\section{A refinement}
\label{refinement}
\setcounter{equation}{0}

        Let us continue with the same notations as in the previous
section.  Observe that
\begin{equation}
        X(\alpha, \eta, C) \subseteq X(\alpha, \delta, C)
\end{equation}
when $\delta \le \eta$.  Put
\begin{equation}
        X(\alpha, C) = \bigcup_{\delta > 0} X(\alpha, \delta, C),
\end{equation}
which is the same as
\begin{equation}
        X(\alpha, C) = \bigcup_{l = 1}^\infty X(\alpha, 1/l, C),
\end{equation}
by monotonicity in $\delta$.  Thus $X(\alpha, C)$ is a Borel set,
since it is the union of a sequence of closed sets.  We would like to
show that
\begin{equation}
\label{mu(A) le C H^alpha(A)}
        \mu(A) \le C \, H^\alpha(A)
\end{equation}
when $A$ is a Borel set such that $A \subseteq X(\alpha, C)$.  Consider
\begin{equation}
        A_l = A \cap X(\alpha, 1/l, C).
\end{equation}
This is a Borel set for each $l$, and
\begin{equation}
        \mu(A_l) \le C \, H^\alpha_{1/l}(A_l) \le C \, H^\alpha(A_l),
\end{equation}
as in the previous section.  Hence
\begin{equation}
        \mu(A_l) \le C \, H^\alpha(A)
\end{equation}
for each $l$, which implies (\ref{mu(A) le C H^alpha(A)}), because $A
= \bigcup_{l = 1}^\infty A_l$ and $A_l \subseteq A_{l + 1}$ for each
$l$.

\section{Hausdorff dimension}
\label{hausdorff dimension}
\setcounter{equation}{0}

        The \emph{Hausdorff dimension} ${\dim}_H A$ of a set $A$ in a
metric space $M$ can normally be defined as the supremum of the set of
$\alpha \ge 0$ such that $H^\alpha(A) = +\infty$, or as the infimum of
the set of $\beta \ge 0$ such that $H^\beta(A) = 0$.  These two
quantities are the same when they are both defined, because
$H^\beta(A) = 0$ when $H^\alpha(A) < + \infty$ and $\beta > \alpha$,
as in Section \ref{hausdorff measures}.  If $A$ has only finitely many
elements, then $H^0(A) < + \infty$ and $H^\alpha(A) = 0$ for every
$\alpha > 0$, and ${\dim}_H A = 0$.  Similarly, ${\dim}_H A = +\infty$
when $H^\alpha(A) = +\infty$ for every $\alpha \ge 0$.

        If $A \subseteq B \subseteq M$, then
\begin{equation}
\label{{dim}_H A le {dim}_H B}
        {\dim}_H A \le {\dim}_H B.
\end{equation}
If $A_1, A_2, \ldots$ is a sequence of subsets of $M$, then
\begin{equation}
        {\dim}_H \Big(\bigcup_{j = 1}^\infty A_j\Big)
          = \sup_{j \ge 1} \, {\dim}_H A_j.
\end{equation}
Indeed, (\ref{{dim}_H A le {dim}_H B}) implies that
\begin{equation}
        {\dim}_H A_i \le {\dim}_H \Big(\bigcup_{j = 1}^\infty A_j\Big)
\end{equation}
for each $i$, and the opposite inequality follows from countable subadditivity.

        Note that ${\dim}_H A = \alpha$ when
\begin{equation}
        0 < H^\alpha(A) < +\infty.
\end{equation}
For example, an interval of positive length in the real line has
Hausdorff dimension $1$, and the Hausdorff dimension of a ball or a
cube in ${\bf R}^n$ is $n$.

\section{Separating sets}
\label{separating sets}
\setcounter{equation}{0}

        Let $A$ be a closed set in ${\bf R}^n$, where the latter is
equipped with the standard Euclidean metric.  If ${\bf R}^n \backslash
A$ is not connected, then
\begin{equation}
        H^{n - 1}(A) > 0.
\end{equation}
To see this, let $p$, $q$ be elements of different connected
components of ${\bf R}^n \backslash A$.  Thus $p \ne q$, and there is
a unique line $L$ passing through them.  Let $P$ be the hyperplane in
${\bf R}^n$ passing through $p$ and perpindicular to $L$, and let
$\pi$ be the orthogonal projection of ${\bf R}^n$ onto $P$.
A key point now is that
\begin{equation}
        \hbox{$\pi(A)$ contains a neighborhood of $p$ in $P$}.
\end{equation}
For if $z \in P \backslash \pi(A)$, then the line $L(z)$ passing
through $z$ and parallel to $L$ is disjoint from $A$.  If $z$ is
sufficiently close to $p$, then $L(z)$ intersects each of the
complementary components of $A$ containing $p$ and $q$, a
contradiction.  It follows that
\begin{equation}
        H^{n - 1}(\pi(A)) > 0.
\end{equation}
The remaining point is that
\begin{equation}
        H^{n - 1}(A) \ge H^{n - 1}(\pi(A)),
\end{equation}
basically because coverings of $A$ can be projected onto coverings of
$\pi(A)$ without increasing the diameters of the sets in the
coverings.

\section{Lipschitz mappings}
\label{lipschitz mappings}
\setcounter{equation}{0}

        Let $(M_1, d_1(x, y))$ and $(M_2, d_2(u, v))$ be metric
spaces.  A mapping $f : M_1 \to M_2$ is said to be \emph{Lipschitz} if
there is a $k \ge 0$ such that
\begin{equation}
\label{d_2(f(x), f(y)) le k d_1(x, y)}
        d_2(f(x), f(y)) \le k \, d_1(x, y)
\end{equation}
for every $x, y \in M_1$.  We may also say that $f$ is $k$-Lipschitz
in this case, to be more precise.  Thus a mapping is $0$-Lipschitz if
and only if it is constant.  If $M_1 = M_2$ with the same metric, then
the identity mapping $f(x) = x$ is $1$-Lipschitz.  Orthogonal
projections of ${\bf R}^n$ onto affine subspaces are also
$1$-Lipschitz.  Another class of examples will be given in the next
section.

        Suppose that $f : M_1 \to M_2$ is $k$-Lipschitz.  If $E
\subseteq M_1$ is bounded, then $f(E)$ is bounded in $M_2$, and
\begin{equation}
        \diam f(E) \le k \, \diam E.
\end{equation}
Using this, one can check that
\begin{equation}
        H^\alpha_{k \delta}(f(A)) \le k^\alpha \, H^\alpha_\delta(A)
\end{equation}
for every $A \subseteq M_1$, $\alpha \ge 0$, and $0 < \delta \le
\infty$, and hence
\begin{equation}
        H^\alpha(f(A)) \le k^\alpha \, H^\alpha(A).
\end{equation}
In particular, ${\dim}_H f(A) \le {\dim}_H A$.

\section{Connected sets}
\label{connected sets}
\setcounter{equation}{0}

        Let $(M, d(x, y))$ be a metric space.  For each $p \in M$,
\begin{equation}
\label{f_p(x) = d(p, x)}
        f_p(x) = d(p, x)
\end{equation}
is $1$-Lipschitz as a mapping from $M$ into ${\bf R}$ with the
standard metric.  This can be verified using the triangle inequality.

        If $A \subseteq M$ is connected, then
\begin{equation}
\label{diam A le H^1(A)}
        \diam A \le H^1(A).
\end{equation}
To see this, we can use the remarks in the previous section to get that
\begin{equation}
        H^1(f_p(A)) \le H^1(A)
\end{equation}
for every $p \in M$.  The connectedness of $A$ and continuity of $f_p$
imply that $f_p(A)$ is basically an interval in ${\bf R}$, whose
$1$-dimensional Hausdorff measure is the same as its length, which is
the same as its diameter, so that
\begin{equation}
        \diam f_p(A) \le H^1(A).
\end{equation}
This implies (\ref{diam A le H^1(A)}), because the diameter of $A$ is
the same as the supremum of the diameter of $f_p(A)$ over $p \in A$.

\section{Bilipschitz enbeddings}
\label{bilipschitz}
\setcounter{equation}{0}

        Let $(M_1, d_1(x, y))$, $(M_2, d_2(u, v))$ be metric spaces
again.  A mapping $f$ from $M_1$ into $M_2$ is said to be
\emph{bilipschitz} if there is a $k \ge 1$ such that
\begin{equation}
\label{k^{-1} d_1(x, y) le d_2(f(x), f(y)) le k d_1(x, y)}
        k^{-1} \, d_1(x, y) \le d_2(f(x), f(y)) \le k \, d_1(x, y)
\end{equation}
for every $x, y \in M_1$.  Equivalently, $f : M_1 \to M_2$ is
$k$-bilipschitz if $f$ is $k$-Lipschitz, $f$ is one-to-one, and the
inverse mapping is $k$-Lipschitz on $f(M_1)$.  In this case,
\begin{equation}
 k^{-\alpha} \, H^\alpha(A) \le H^\alpha(f(A)) \le k^\alpha \, H^\alpha(A)
\end{equation}
for every $A \subseteq M_1$ and $\alpha \ge 0$, and in particular
\begin{equation}
        {\dim}_H f(A) = {\dim}_H A
\end{equation}
for every $A \subseteq M$.

        By contrast, Hausdorff dimensions can be changed by ordinary
homeomorphisms.  For instance, there are topological Cantor sets in
the real line with any Hausdorff dimension $\alpha$, $0 \le \alpha \le
1$.  There are even Cantor sets with positive Lebesgue measure, and
hence positive one-dimensional Hausdorff measure.

\section{Ultrametric spaces}
\label{ultrametrics}
\setcounter{equation}{0}

        Let $(M, d(x, y))$ be a metric space.  We say that $d(x, y)$
is an \emph{ultrametric} on $M$ if
\begin{equation}
\label{d(x, z) le max(d(x, y), d(y, z))}
        d(x, z) \le \max(d(x, y), d(y, z))
\end{equation}
for every $x, y, z \in M$.  Of course, this is stronger than the usual
triangle inequality.

        A basic class of examples of ultrametric spaces is given by
sequence spaces.  Let $L$ be a nonempty set, and let $X$ be the set of
sequences $x = \{x_j\}_{j = 1}^\infty$ with $x_j \in L$ for each $j$.
Also let $\rho$ be a positive real number, $\rho < 1$.  Let $d_\rho(x,
y)$ be defined for $x, y \in X$ by $d_\rho(x, y) = 0$ when $x = y$, and
\begin{equation}
        d_\rho(x, y) = \rho^l
\end{equation}
where $l$ is the largest nonnegative integer such that $x_j = y_j$ for
$j \le l$ otherwise.  If $x_1 \ne y_1$, then $l = 0$, and $d_\rho(x,
y) = 1$.  It is not difficult to check that $d_\rho(x, y)$ satisfies
the ultrametric version of the triangle inequality.  Basically, this
amounts to the statement that if $x, y, z \in X$, and $x_j = y_j$,
$y_j = z_j$ when $j \le l$, then $x_j = z_j$ when $j \le l$.  The same
argument works for $\rho = 1$, in which case $d_\rho(x, y)$ reduces to
the discrete metric on $X$.  The topology on $X$ determined by
$d_\rho(x, y)$ is the same as the product topology when $\rho < 1$,
where $X$ is considered as the product of infinitely many copies of
$L$, and $L$ is equipped with the discrete topology.  In particular,
$X$ is compact when $L$ has only finitely many elements.

        If $L$ is a finite set with at least two elements, then $X$ is
a topological Cantor set, which is to say that it is homeomorphic to
the usual middle-thirds Cantor set.  If $L$ has exactly two elements
and $\rho = 1/3$, then $X$ is bilipschitz equivalent to the
middle-thirds Cantor set.

\section{Closed balls}
\label{closed balls}
\setcounter{equation}{0}

        Let $(M, d(x, y))$ be a metric space.  If $p \in M$ and $r \ge 0$,
then the closed ball in $M$ with center $p$ and radius $r$ is defined by
\begin{equation}
        \overline{B}(p, r) = \{x \in M : d(p, x) \le r\}.
\end{equation}
This is a closed and bounded set in $M$, with
\begin{equation}
        \diam \overline{B}(p, r) \le 2 \, r.
\end{equation}

        Suppose now that $d(x, y)$ is an ultrametric on $M$.  In this
case, $\overline{B}(p, r)$ is actually an open set in $M$ too, because
\begin{equation}
\label{B(q, r) subseteq overline{B}(q, r) subseteq overline{B}(p, r)}
        B(q, r) \subseteq \overline{B}(q, r) \subseteq \overline{B}(p, r)
\end{equation}
for every $q \in \overline{B}(p, r)$.  Similarly, open balls in $M$
are closed sets.  Moreover,
\begin{equation}
        \diam \overline{B}(p, r) \le r.
\end{equation}

        Let us say that $\mathcal{C}$ is a \emph{cell} in an
ultrametric space $M$ if it is a closed ball, i.e., if it can be
represented as $\overline{B}(p, r)$ for some $p \in M$ and $r \ge 0$.
In some circumstances, one may prefer to use only closed balls of
positive radius.  If $E$ is a bounded set in $M$, $p \in E$, and
\begin{equation}
        \mathcal{C} = \overline{B}(p, \diam E),
\end{equation}
then $E \subseteq \mathcal{C}$ and
\begin{equation}
        \diam \mathcal{C} = \diam E.
\end{equation}
This permits one to use coverings by cells in the definition of
Hausdorff measure in an ultrametric space.

        If $\mathcal{C}$ is a cell in an ultrametric space $M$, then
\begin{equation}
        \mathcal{C} = \overline{B}(p, \diam \mathcal{C})
\end{equation}
for every $p \in \mathcal{C}$.  If $\mathcal{C}_1$, $\mathcal{C}_2$
are cells in $M$ such that
\begin{equation}
        \mathcal{C}_1 \cap \mathcal{C}_2 \ne \emptyset,
\end{equation}
then
\begin{equation}
        \mathcal{C}_1 \subseteq \mathcal{C}_2
            \quad\hbox{or}\quad
        \mathcal{C}_2 \subseteq \mathcal{C}_1.
\end{equation}
This follows by representing $\mathcal{C}_1$, $\mathcal{C}_2$ as
closed balls centered at the same point $p$.

\section{Sequence spaces}
\label{sequence spaces}
\setcounter{equation}{0}

        Let $L$ be a finite set with $n \ge 2$ elements, let $X$ be
the space of sequences with entries in $L$, and let $d_\rho(x, y)$ be
the ultrametric on $X$ associated to some $\rho \in (0, 1)$ as in
Section \ref{ultrametrics}.  If $\alpha$ is the positive real number
defined by
\begin{equation}
        n \, \rho^\alpha = 1,
\end{equation}
then
\begin{equation}
        H^\alpha(X) = H^\alpha_{con}(X) = 1.
\end{equation}
To see this, observe first that
\begin{equation}
        H^\alpha(X) \le 1,
\end{equation}
since $X$ can be covered by $n^l$ cells of diameter $\rho^l$ for each $l$.
Similarly, every cell in $X$ of diameter $\rho^j$ can be covered by
$n^{l - j}$ cells of diameter $\rho^l$ when $j \le l$.  Using this
fact, one can replace any covering of $X$ by finitely many cells
$\mathcal{C}_1, \ldots, \mathcal{C}_k$ with a covering by subcells
$\mathcal{C}'_1, \ldots, \mathcal{C}'_r$ in such a way that
\begin{equation}
        \diam \mathcal{C}'_i = \rho^l
\end{equation}
for some $l$ and all $i = 1, \ldots, r$, and
\begin{equation}
        \sum_{h = 1}^k (\diam \mathcal{C}_h)^\alpha
         = \sum_{i = 1}^r (\diam \mathcal{C}'_i)^\alpha.
\end{equation}
In order for $X$ to be covered by a collection of cells of the same
diameter $\rho^l$, it is necessary for all $n^l$ of these cells to be
used.  This implies that
\begin{equation}
        \sum_{h = 1}^k (\diam \mathcal{C}_h)^\alpha
         = \sum_{i = 1}^r (\diam \mathcal{C}'_i)^\alpha \ge 1,
\end{equation}
and hence
\begin{equation}
        H^\alpha(X) \ge H^\alpha_{con}(X) \ge 1.
\end{equation}

\section{Snowflake metrics}
\label{snowflakes}
\setcounter{equation}{0}

        If $a, b \ge 0$ and $0 < t \le 1$, then
\begin{equation}
\label{(a + b)^t le a^t + b^t}
        (a + b)^t \le a^t + b^t.
\end{equation}
Indeed,
\begin{equation}
        \max(a, b) \le (a^t + b^t)^{1/t},
\end{equation}
which implies that
\begin{eqnarray}
        a + b & \le & \max(a, b)^{1 - t} \, (a^t + b^t) \\
 & \le & (a^t + b^t)^{(1 - t)/t} \, (a^t + b^t) = (a^t + b^t)^{1/t}. \nonumber
\end{eqnarray}

        Let $(M, d(x, y))$ be a metric space.  It follows from
(\ref{(a + b)^t le a^t + b^t}) that $d(x, y)^t$ is also a metric on
$M$ for each $t \in (0, 1)$.  If $d(x, y)$ is an ultrametric on $M$,
then it is easy to see that $d(x, y)^t$ is an ultrametric on $M$ for
every $t > 0$.  It is easy to keep track of the affect of this change
in the metric on diameters of subsets of $M$ and so on as well.  In
particular, the $\alpha$-dimensional Hausdorff measure on $M$ with
respect to the initial metric $d(x, y)$ is the same as Hausdorff
measure of dimension $\alpha / t$ with respect to the new metric $d(x,
y)^t$ for each $\alpha \ge 0$.

\section{Topological dimension $0$}
\label{topological dimension 0}
\setcounter{equation}{0}

        As in \cite{h-w}, a separable metric space $M$ is said to have
\emph{topological dimension $0$} if for every $p \in M$ and $r > 0$
there is an open set $U$ in $M$ such that $p \in U \subseteq B(p, r)$
and $\partial U = \emptyset$.  Equivalently, $M$ has topological
dimension $0$ if for every $p \in M$ and open set $V$ in $M$ with $p
\in V$ there is an open set $U$ in $M$ such that $p \in U \subseteq V$
and $\partial U = \emptyset$, so that this condition depends only on
the topology on $M$.  If the metric $d(x, y)$ on $M$ is an
ultrametric, then every open ball in $M$ is a closed set, and $M$ has
topological dimension $0$.  For any metric $d(x, y)$ on $M$,
\begin{equation}
        H^1(M) = 0
\end{equation}
implies that $M$ has topological dimension $0$.  Remember that
\begin{equation}
        f_p(x) = d(p, x)
\end{equation}
is a $1$-Lipschitz function for each $p \in M$, and hence
\begin{equation}
        H^1(f_p(M)) \le H^1(M) = 0.
\end{equation}
Here $H^1(f_p(M))$ is the $1$-dimensional Hausdorff measure of
$f_p(M)$ as a set in the real line, with the standard metric.  Thus,
for each $p \in M$,
\begin{equation}
        \{x \in M : d(p, x) = r\} = \emptyset
\end{equation}
for almost every $r$ with respect to Lebesgue measure on ${\bf R}$,
which implies that $M$ has topological dimension $0$.

\section{Distance functions}
\label{distance functions}
\setcounter{equation}{0}

        Let $(M, d(x, y))$ be a metric space.  If $A \subseteq M$, $A
\ne \emptyset$, and $x \in M$, then put
\begin{equation}
        \dist(x, A) = \inf \{d(x, a) : a \in A\}.
\end{equation}
Thus $\dist(x, A) = 0$ if and only if $x \in \overline{A}$.  Also, the
distance from any $x \in M$ to $A$ is the same as the distance from
$x$ to the closure of $A$.  It is easy to check that
\begin{equation}
        \dist(x, A) \le \dist(y, A) + d(x, y)
\end{equation}
for every $x, y \in M$, so that $\dist(x, A)$ is $1$-Lipschitz in $x$.
If $d(x, y)$ is an ultrametric on $M$, then $\dist(x, M)$ satisfies
the ``ultra-Lipschitz'' property
\begin{equation}
        \dist(x, A) \le \max(\dist(y, A), d(x, y)).
\end{equation}
In particular,
\begin{equation}
        \dist(x, A) = \dist(y, A)
\end{equation}
when $d(x, y) < \dist(x, A)$.

        The continuity of $\dist(x, A)$ implies that
\begin{equation}
        U_r = \{x \in M : \dist(x, A) < r\}
\end{equation}
is an open set in $M$ for each $r > 0$, and
\begin{equation}
        L_r = \{x \in M : \dist(x, A) \le r\}
\end{equation}
is a closed set.  If $d(x, y)$ is an ultrametric, then $U_r$ is a
closed set for each $r > 0$, and $L_r$ is an open set.  For any metric
$d(x, y)$, $H^1(M) = 0$ implies that
\begin{equation}
        \{x \in M : \dist(x, A) = r\} = \emptyset
\end{equation}
for almost every $r > 0$, so that $U_r = L_r$ for almost every $r$ in
this case.

\section{Disjoint closed sets}
\label{disjoint closed sets}
\setcounter{equation}{0}

        Let $(M, d(x, y))$ be a metric space, and let $A$, $B$ be
disjoint nonempty closed subsets of $M$.  Consider the function
\begin{equation}
        \phi(x) = \frac{\dist(x, A)}{\dist(x, A) + \dist(x, B)}.
\end{equation}
As usual, this is a continuous real-valued function on $M$ such that
$\phi(x) = 0$ when $x \in A$, $\phi(x) = 1$ when $x \in B$, and $0 \le
\phi(x) \le 1$ for every $x \in M$.  For each $\epsilon > 0$, $\phi$
is also Lipschitz on the set
\begin{equation}
 \Lambda_\epsilon = \{x \in M : \dist(x, A) + \dist(x, B) \ge \epsilon\}.
\end{equation}
In particular, $\phi$ is locally Lipschitz, in the sense that $\phi$
is Lipschitz on a neighborhood of any point in $M$.  If the distance
between $A$ and $B$ is positive, then $\phi$ is Lipschitz on all of
$M$.  Of course, this holds automatically when $A$ or $B$ is compact.
If $d(x, y)$ is an ultrametric on $M$, then $\dist(x, A)$ is locally
constant on $M \backslash A$, $\dist(x, B)$ is locally constant on $M
\backslash B$, and hence $\phi$ is locally constant on $M \backslash
(A \cup B)$.

        For any metric $d(x, y)$ on $M$, if $H^1(M) = 0$, then
\begin{equation}
        H^1(\phi(\Lambda_\epsilon)) = 0,
\end{equation}
for each $\epsilon > 0$, since $\phi$ is Lipschitz on $\Lambda_\epsilon$.
This implies that
\begin{equation}
        H^1(\phi(M)) = 0,
\end{equation}
since $M = \bigcup_{n = 1}^\infty \Lambda_{1/n}$.  Thus $\phi^{-1}(r)
= \emptyset$ for almost every $r \in (0, 1)$.

\section{Topological dimension $0$, 2}
\label{topological dimension 0, 2}
\setcounter{equation}{0}

        Let $M$ be a separable metric space.  It is easy to see that
$M$ has topological dimension $0$ if and only if for each $p \in M$
and closed set $B$ in $M$ with $p \not\in B$ there is an open set $U
\subseteq M$ such that
\begin{equation}
        p \in U, \, U \cap B = \emptyset, \hbox{ and } \partial U = \emptyset.
\end{equation}
Let us say that $M$ has property P if for every pair of distinct
elements $p$, $q$ of $M$ there is an open set $U$ in $M$ such that
\begin{equation}
        p \in U, \, q \not\in U, \hbox{ and } \partial U = \emptyset.
\end{equation}
Similarly, we say that $M$ has property Q if for every pair $A$, $B$
of disjoint closed subsets of $M$ there is an open set $U$ in $M$ such that
\begin{equation}
 A \subseteq U, \, U \cap B = \emptyset, \hbox{ and } \partial U = \emptyset.
\end{equation}
Thus property Q implies that $M$ has topological dimension $0$, which
implies property P, which implies that $M$ is totally disconnected in
the sense that it has no connected subsets with at least two elements.
Note that properties P and Q are automatically symmetric in $p$, $q$
and $A$, $B$, respectively, since their roles in the definitions can
be interchanged.  The conclusions of properties P and Q are also
symmetric in $p$, $q$ and $A$, $B$, because $V = M \backslash U$ is an
open set in $M$ with $\partial V = \emptyset$ when $U$ is an open set
in $M$ with $\partial U = \emptyset$.

        Suppose that $M$ satisfies property P, $p \in M$, $B$ is a
compact set in $M$, and $p \not\in B$.  For each $q \in B$, let $V(q)$
be an open set in $M$ such that $q \in V(q)$, $p \not\in V(q)$, and
$\partial V(q) = \emptyset$.  By compactness, there are finitely many
elements $q_1, \ldots, q_l$ of $B$ such that
\begin{equation}
        B \subseteq V(q_1) \cup V(q_2) \cup \cdots \cup V(q_l).
\end{equation}
Hence $V = \bigcup_{i = 1}^l V(q_i)$ is an open set in $M$ such that
$B \subseteq V$, $p \not\in V$, and $\partial V = \emptyset$.
Equivalently, $U = M \backslash V$ is an open set in $M$ such that $p
\in U$, $B \subseteq M \backslash U$, and $\partial U = \emptyset$.
If $A$, $B$ are disjoint compact subsets of $M$, then one can apply
this to each $p \in A$ to get an open set $U(p)$ in $M$ such that $p
\in U(p)$, $B \subseteq M \backslash U(p)$, and $\partial U(p) =
\emptyset$.  By compactness, there are finitely many elements $p_1,
\ldots, p_r$ of $A$ such that
\begin{equation}
        A \subseteq U(p_1) \cup U(p_2) \cup \cdots \cup U(p_r).
\end{equation}
Thus $W = \bigcup_{j = 1}^r U(p_j)$ is an open set, $A \subseteq W$,
$B \subseteq M \backslash W$, and $\partial W = \emptyset$.  This
shows that property P implies the analogues of the reformulation of
topological dimension $0$ in the previous paragraph and property Q for
compact subsets of $M$ instead of closed subsets of $M$.  In
particular, property P implies that $M$ has topological dimension $0$
and satisfies property Q when $M$ is compact.

        The same type of covering argument shows that $M$ satisfies
the analogue of property Q for compact sets $A$ in $M$ when $M$ has
topological dimension $0$.  An important theorem states that any separable
metric space of topological dimension $0$ satisfies property Q.  The proof
of this will be given in Section \ref{topological dimension 0, 3}.

        Let $A$, $B$ be disjoint nonempty closed subsets of a metric
space $M$, and let $\phi(x)$ be as in the previous section.  Thus
\begin{equation}
        U_r = \{x \in M : \phi(x) < r\}
\end{equation}
is an open set in $M$ such that
\begin{equation}
        A \subseteq U_r \hbox{ and } \overline{U_r} \cap B = \emptyset
\end{equation}
for every $r \in (0, 1)$.  If the metric $d(x, y)$ on $M$ is an
ultrametric, then $\partial U = \emptyset$ for each $r \in (0, 1)$,
because $\phi$ is locally constant on $M \backslash (A \cup B)$.
For any metric on $M$,
\begin{equation}
        \partial U_r \subseteq \{x \in M : \phi(x) = r\}.
\end{equation}
If $H^1(M) = 0$, then it follows that $\partial U_r = \emptyset$ for
almost every $r \in (0, 1)$, so that $M$ has property Q.

\section{Strong normality}
\label{strong normality}
\setcounter{equation}{0}

        Remember that two subsets $A$, $B$ of a topological space $X$
are said to be \emph{separated} if
\begin{equation}
        \overline{A} \cap B = A \cap \overline{B} = \emptyset.
\end{equation}
A strong version of normality asks that for any pair of separated
subsets $A$, $B$ of $X$ there are disjoint open subsets $U$, $V$ of
$X$ such that $A \subseteq U$ and $B \subseteq V$.  This implies that
every $Y \subseteq X$ has the same property with respect to the
induced topology, because $A, B \subseteq Y$ are separated relative to
$Y$ if and only if they are separated relative to $X$.

        If the topology on $X$ is defined by a metric, then $X$
satisfies this stronger version of normality.  To see this, let $A$,
$B$ be separated subsets of $X$.  By hypothesis, each $a \in A$ is not
an element of the closure of $B$, and so there is an $r(a) > 0$ such that
\begin{equation}
        B(a, r(a)) \cap B = \emptyset.
\end{equation}
Similarly, for every $b \in B$ there is a $t(b) > 0$ such that
\begin{equation}
        B(b, t(b)) \cap A = \emptyset.
\end{equation}
Consider
\begin{equation}
\label{U = bigcup_{a in A} B(a, r(a)/2), V = bigcup_{b in B} B(b, t(b)/2)}
 U = \bigcup_{a \in A} B(a, r(a)/2), \ V = \bigcup_{b \in B} B(b, t(b)/2).
\end{equation}
Thus $U$ and $V$ are open subsets of $X$ such that $A \subseteq U$ and
$B \subseteq V$.  If $U \cap V \ne \emptyset$, then there are $a \in
A$ and $b \in B$ such that
\begin{equation}
        B(a, r(a)/2) \cap B(b, t(b)/2) \ne \emptyset.
\end{equation}
This would imply that the distance from $a$ to $b$ is less than the
average of $r(a)$ and $t(b)$, by the triangle inequality.  However,
the distance from $a$ to $b$ is greater than or equal to both $r(a)$
and $t(b)$ by construction, and so we may conclude that $U$ and $V$
are disjoint.  Note that one can use balls of radius $r(a)$ and $t(b)$
in (\ref{U = bigcup_{a in A} B(a, r(a)/2), V = bigcup_{b in B} B(b,
t(b)/2)}) when the metric on $X$ is an ultrametric.

        A classical theorem asserts that a regular topological space
$X$ with a countable base for its topology satisfies this strong
version of normality.  Indeed, let $A$, $B$ be separated subsets of
$X$ again.  Regularity implies that each $a \in A$ has a neighborhood
$U(a)$ such that
\begin{equation}
        \overline{U(a)} \cap B = \emptyset.
\end{equation}
Similarly, each $b \in B$ has a neighborhood $V(b)$ such that
\begin{equation}
        \overline{V(b)} \cap A = \emptyset.
\end{equation}
Because $X$ has a countable base for its topology, there are sequences
$\{U_i\}_{i = 1}^\infty$ and $\{V_j\}_{j = 1}^\infty$ of open subsets
of $X$ such that
\begin{equation}
        A \subseteq \bigcup_{i = 1}^\infty U_i,
           \ B \subseteq \bigcup_{j = 1}^\infty V_j
\end{equation}
and
\begin{equation}
        \overline{U_i} \cap B = \overline{V_j} \cap A = \emptyset
\end{equation}
for every $i, j$.  For each $l, n \ge 1$, put
\begin{equation}
\widetilde{U}_l = U_l \backslash \Big(\bigcup_{j = 1}^l \overline{V_j}\Big), \ 
 \widetilde{V}_n = V_n \backslash \Big(\bigcup_{i = 1}^n \overline{U_i}\Big)
\end{equation}
Thus $\widetilde{U}_l$, $\widetilde{V}_n$ are open sets for every $l,
n \ge 1$, and
\begin{equation}
        \widetilde{U}_l \cap \widetilde{V}_n = \emptyset.
\end{equation}
This implies that 
\begin{equation}
        U = \bigcup_{l = 1}^\infty \widetilde{U}_l,
         \ V = \bigcup_{n = 1}^\infty \widetilde{V}_n
\end{equation}
are disjoint open sets, and it is easy to check that $A \subseteq U$,
$B \subseteq V$, as desired.  Of course, Urysohn's famous metrization
theorem implies that a regular topological space with a countable base
for its topology is metrizable.

\section{Topological dimension $0$, 3}
\label{topological dimension 0, 3}
\setcounter{equation}{0}

        Let $M$ be a separable metric space of topological dimension
$0$, let $A$, $B$ be disjoint closed subsets of $M$, and let us show
that there is an open and closed set $U$ in $M$ such that $A \subseteq
U$ and $B \subseteq M \backslash U$, as in Section \ref{topological
dimension 0, 2}.  For each $p \in M$, there is an open and closed set
$W(p)$ in $M$ such that $p \in W(p)$ and
\begin{equation}
        W(p) \cap A = \emptyset \hbox{ or } W(p) \cap B = \emptyset.
\end{equation}
Because $M$ has a countable base for its topology, it follows that
there is a sequence $W_1, W_2, \ldots$ of open and closed subsets of
$M$ such that
\begin{equation}
        \bigcup_{i = 1}^\infty W_i = M,
\end{equation}
and
\begin{equation}
        W_i \cap A = \emptyset \hbox{ or } W_i \cap B = \emptyset
\end{equation}
for each $i$.  Put $\widetilde{W}_1 = W_1$ and
\begin{equation}
        \widetilde{W}_l = W_l \backslash \Big(\bigcup_{i = 1}^{l - 1} W_i\Big)
\end{equation}
when $l \ge 2$.  Thus $\widetilde{W}_l \subseteq W_l$ and $\widetilde{W}_l$
is both open and closed for each $l$, and
\begin{equation}
 \bigcup_{l = 1}^\infty \widetilde{W}_l = \bigcup_{i = 1}^\infty W_i = M.
\end{equation}
Consider
\begin{equation}
 U = \bigcup \{\widetilde{W}_l : \widetilde{W}_l \cap A \ne \emptyset\}
\end{equation}
and
\begin{equation}
 V = \bigcup \{\widetilde{W}_l : \widetilde{W}_l \cap A = \emptyset\}
   = M \backslash U.
\end{equation}
These are obviously open sets, since they are unions of open sets.
Hence they are closed sets as well, because they are complements of
each other.  Each $\widetilde{W}_l$ can intersect at most one of $A$
and $B$, which implies that $A \subseteq U$ and $B \subseteq V$, as
desired.

\section{Subsets}
\label{subsets}
\setcounter{equation}{0}

        Let $M$ be a separable metric space.  A set $X \subseteq M$ is
considered to have topological dimension $0$ if it has topological
dimension $0$ as a metric space itself, with respect to the
restriction of the metric on $M$ to $X$.  Actually, it is sometimes
convenient to let the topological dimension of the empty set be $-1$,
but this is not important for the moment.  It is easy to see that $X
\subseteq M$ has topological dimension $0$ when $M$ has topological
dimension $0$.  Hence $X \subseteq Y \subseteq M$ has topological
dimension $0$ when $Y$ has topological dimension $0$.

        Consider the case of the real line with the standard metric.
It is easy to see that $X \subseteq {\bf R}$ has topological dimension
$0$ if and only if $X$ is totally disconnected, which is the same as
saying that $X$ has at least two elements but does not contain any
interval of positive length in this case.  For example, the set ${\bf
Q}$ of rational numbers has topological dimension $0$, as does the set
${\bf R} \backslash {\bf Q}$ of irrational numbers.  This shows that
the union of two sets with topological dimension $0$ may not have
topological dimension $0$, since the real line is connected and
therefore does not have topological dimension $0$.

        Let $M$ be a separable metric space again, and suppose that $X
\subseteq M$ has topological dimension $0$ and $M = X \cup \{p\}$ for
some $p \in M \backslash X$.  Thus $X$ is an open set in $M$, and we
would like to show that $M$ also has topological dimension $0$ in
these circumstances.  It is easy to see that every element of $X$ has
arbitrarily small open neighborhoods with empty boundary in $M$,
because of the corresponding property of $X$.  The main point is to
show that the analogous statement holds at $p$.  Let $r > 0$ be given,
and consider
\begin{equation}
        A = \overline{B}(p, r/2), \ B = M \backslash B(p, r).
\end{equation}
These are disjoint closed subsets of $M$, and $B \subseteq X$.  Hence
$A \backslash \{p\}$, $B$ are disjoint relatively closed subsets of
$X$.  The theorem in the previous section may be applied in $X$ to get
disjoint open sets $U$, $V$ in $X$ and therefore in $M$ such that $A
\backslash \{p\} \subseteq U$, $B \subseteq V$, and $U \cup V = X$.
If $U_1 = U \cup \{p\}$, then $U_1$ is an open set in $M$, because
$B(p, r/2) \subseteq U_1$ by construction.  Moreover, $U_1 \subseteq
B(p, r)$ because $U_1 \cap B = \emptyset$, and $U_1$ is closed in $M$,
since $M \backslash U_1 = V$ is open.

        It follows that the union of a set with topological dimension
$0$ and a finite set also has topological dimension $0$.  The
analogous statement for countable sets does not work, as in the
example of the sets of rational and irrational numbers in the real
line.

\section{Closed subsets}
\label{closed subsets}
\setcounter{equation}{0}

        Let $M$ be a separable metric space, and let $X_1 \subseteq M$
be a closed set with topological dimension $0$.  If $A$, $B$ are
disjoint closed subsets of $M$, then $A \cap X_1$, $B \cap X_1$ are
disjoint closed subsets of $X_1$.  The theorem in Section
\ref{topological dimension 0, 3} implies that there are disjoint
closed subsets $\widetilde{A}$, $\widetilde{B}$ of $X_1$ such that
\begin{equation}
 A \cap X_1 \subseteq \widetilde{A}, \ B \cap X_1 \subseteq \widetilde{B},
                       \hbox{ and } \widetilde{A} \cup \widetilde{B} = X_1.
\end{equation}
Thus $A_1 = A \cup \widetilde{A}$ and $B_1 = B \cup \widetilde{B}$ are
disjoint closed subsets of $M$ such that 
\begin{equation}
\label{A_1, B_1}
        A \subseteq A_1, \ B \subseteq B_1,
           \hbox{ and } X_1 \subseteq A_1 \cup B_1.
\end{equation}

        If $X_2$ is another closed set in $M$ with topological
dimension $0$, then we can apply the same procedure a second time to
get disjoint closed subsets $A_2$, $B_2$ of $M$ such that
\begin{equation}
        A \subseteq A_1 \subseteq A_2, \ B \subseteq B_1 \subseteq B_2,
           \hbox{ and } X_1 \cup X_2 \subseteq A_2 \cup B_2.
\end{equation}
If $M = X_1 \cup X_2$, then $M = A_2 \cup B_2$, so that $A_2$, $B_2$
are open sets as well.  This shows that the union of two closed sets
with topological dimension $0$ also has topological dimension $0$.

        Suppose now that $X_1, X_2, \ldots$ is a sequence of closed
subsets of a separable metric space $M$, where each $X_i$ has
topological dimension $0$.  If $M = \bigcup_{i = 1}^\infty X_i$, then
a fundamental theorem states that $M$ also has topological dimension
$0$.  To prove this, it is not quite sufficient to repeat the process
and take the union of the resulting $A_i$'s and $B_i$'s, because it is
not clear that these will be closed sets.  Instead, one can apply normality
after the first step to get open subsets $U_1$, $V_1$ of $M$ such that
\begin{equation}
\label{U_1, V_1}
        A_1 \subseteq U_1, \ B_1 \subseteq V_1, \hbox{ and }
         \overline{U_1} \cap \overline{V_1} = \emptyset.
\end{equation}
Proceeding as before, we get disjoint closed subsets $\widehat{A}_2$,
$\widehat{B}_2$ of $M$ such that
\begin{equation}
        \overline{U_1} \subseteq \widehat{A}_2,
         \ \overline{V_1} \subseteq \widehat{B}_2,
          \hbox{ and } X_2 \subseteq \widehat{A}_2 \cup \widehat{B}_2.
\end{equation}
Of course,
\begin{equation}
        X_1 \subseteq A_1 \cup B_1 \subseteq U_1 \cup V_1 
             \subseteq \widehat{A}_2 \cup \widehat{B}_2.
\end{equation}
Repeating the process, we get increasing sequences $\{U_i\}_{i =
1}^\infty$, $\{V_i\}_{i = 1}^\infty$ of open subsets of $M$, and we
put $U = \bigcup_{i = 1}^\infty U_i$, $V = \bigcup_{i = 1}^\infty
V_i$.  By construction,
\begin{equation}
        A \subseteq U, \ B \subseteq V, \ U \cap V = \emptyset,
           \hbox{ and } U \cup V = M.
\end{equation}

        Suppose instead that $M = X_1 \cup X_2$, where $X_1$, $X_2$
have topological dimension $0$ and $X_1$ is closed.  We may as well
take $X_2 = M \backslash X_1$, because a subset of a set with
topological dimension $0$ also has topological dimension $0$, so that
$X_2$ is an open set in $M$.  Any open set in a metric space is a
countable union of closed sets, which implies that $X_2$ is a
countable union of closed sets of topological dimension $0$.  Hence
$M$ is a countable union of closed sets of topological dimension $0$,
and therefore has topological dimension $0$ as well.  Alternatively,
we can start with a pair of disjoint closed sets $A$, $B$ in $M$, and
apply the earlier argument to get disjoint closed sets $A_1$, $B_1$ in
$M$ satisfying (\ref{A_1, B_1}).  By normality, there are open subsets
$U_1$, $V_1$ of $M$ that satisfy (\ref{U_1, V_1}), so that
$\overline{U_1} \cap X_2$, $\overline{V_1} \cap X_2$ are disjoint
relatively closed sets in $X_2$.  The theorem in Section
\ref{topological dimension 0, 3} implies that there are disjoint open
sets $W$, $Z$ in $X_2$ such that
\begin{equation}
 \overline{U_1} \cap X_2 \subseteq W, \ \overline{V_1} \cap X_2 \subseteq Z,
  \hbox{ and } W \cup Z = X_2.
\end{equation}
Thus $U_1 \cup W$, $V_1 \cup Z$ are disjoint open subsets of $M$ such that
\begin{equation}
        A \subseteq U_1 \cup W, \ B \subseteq V_1 \cup Z,
           \hbox{ and } (U_1 \cup W) \cup (V_1 \cup Z) = M.
\end{equation}

\section{Extrinsic conditions}
\label{extrinsic conditions}
\setcounter{equation}{0}

        Let $M$ be a separable metric space, suppose that $X \subseteq
M$ has topological dimension $0$, and let $A$, $B$ be disjoint closed
subsets of $M$.  We would like to show that there is an open set $W
\subseteq M$ such that
\begin{equation}
\label{A, B, W}
        A \subseteq W, \ \overline{W} \cap B = \emptyset,
           \hbox{ and } \partial W \cap X = \emptyset.
\end{equation}
If $X$ is a closed set in $M$, then one can get disjoint closed
subsets $A_1$, $B_1$ of $M$ as in (\ref{A_1, B_1}), and use normality
to get an open set $W \subseteq M$ such that
\begin{equation}
        A_1 \subseteq W \hbox{ and } \overline{W} \cap B_1 = \emptyset,
\end{equation}
and which therefore satisfies (\ref{A, B, W}).  If $H^1(X) = 0$, then
one can get an open set $W$ that satisfies (\ref{A, B, W}) using the
function $\phi$ in Section \ref{disjoint closed sets}, as in Section
\ref{topological dimension 0, 2}.

        By normality, there are open subsets $U$, $V$ of $M$ such that
\begin{equation}
        A \subseteq U, \ B \subseteq V, \hbox{ and }
                 \overline{U} \cap \overline{V} = \emptyset.
\end{equation}
Thus $\overline{U} \cap X$, $\overline{V} \cap X$ are disjoint
relatively closed subsets of $X$.  The theorem in Section
\ref{topological dimension 0, 3} implies that there are disjoint
relatively open and closed subsets $C$, $D$ of $X$ such that
\begin{equation}
        \overline{U} \cap X \subseteq C, \ \overline{V} \cap X \subseteq D,
         \hbox{ and } C \cup D = X.
\end{equation}
In particular,
\begin{equation}
\label{overline{C} cap D = C cap overline{D} = emptyset}
        \overline{C} \cap D = C \cap \overline{D} = \emptyset
\end{equation}
and
\begin{equation}
        \overline{U} \cap D = C \cap \overline{V} = \emptyset.
\end{equation}
The latter implies that $D \subseteq M \backslash U$, $C \subseteq M
\backslash V$, and hence
\begin{equation}
\overline{D} \subseteq M \backslash U, \ \overline{C} \subseteq M \backslash V,
\end{equation}
since $U$, $V$ are open subsets of $M$.  Equivalently,
\begin{equation}
        U \cap \overline{D} = \overline{C} \cap V = \emptyset,
\end{equation}
which implies that
\begin{equation}
\label{A cap overline{D} = overline{C} cap B = emptyset}
        A \cap \overline{D} = \overline{C} \cap B = \emptyset.
\end{equation}

        Note that
\begin{equation}
        (A \cup C) \cap (B \cup D) = \emptyset.
\end{equation}
Using (\ref{overline{C} cap D = C cap overline{D} = emptyset}) and
(\ref{A cap overline{D} = overline{C} cap B = emptyset}), one can
check that $A \cup C$ and $B \cup D$ are separated, i.e.,
\begin{equation}
        (\overline{A \cup C}) \cap (B \cup D) = \emptyset, \ 
        (A \cup C) \cap (\overline{B \cup D}) = \emptyset.
\end{equation}
The strong version of normality for metric spaces implies that
there is an open set $W \subseteq M$ such that
\begin{equation}
        A \cup C \subseteq W, \ \overline{W} \cap (B \cup D) = \emptyset.
\end{equation}
This implies (\ref{A, B, W}), because $C \cup D = X$.

\section{Level sets}
\label{level sets}
\setcounter{equation}{0}

        Let $(M, d(x, y))$ be a metric space, and let $f$ be a
real-valued $k$-Lipschitz function on $M$.  Suppose that $A \subseteq
M$, $\alpha \ge 1$, and $H^\alpha_{con}(A) < + \infty$.  Let $\{E_i\}_i$
be a collection of finitely or countably many subsets of $M$ such that
\begin{equation}
        A \subseteq \bigcup_i E_i
\end{equation}
and
\begin{equation}
        \sum_i (\diam E_i)^\alpha < + \infty.
\end{equation}
For each $i$, let $\chi_i(r)$ be the characteristic function of
$\overline{f(E_i)}$ on the real line, equal to $1$ when $r \in
\overline{f(E_i)}$ and $0$ when $r \in {\bf R} \backslash
\overline{f(E_i)}$.  Consider
\begin{equation}
        h(r) = \sum_i (\diam E_i)^{\alpha - 1} \, \chi_i(r).
\end{equation}
This is a measurable function on ${\bf R}$, since each $\chi_i$ is
measurable.  If $|B|$ denotes the Lebesgue measure of $B \subseteq
{\bf R}$, then
\begin{equation}
        |\overline{f(E_i)}| \le \diam \overline{f(E_i)} \le k \, \diam E_i
\end{equation}
for each $i$, and hence
\begin{equation}
        \int_{\bf R} h(r) \, dr
          = \sum_i (\diam E_i)^{\alpha - 1} \, |\overline{f(E_i)}|
          \le k \sum_i (\diam E_i)^\alpha.
\end{equation}
For each $r \in {\bf R}$,
\begin{equation}
        A \cap f^{-1}(r) \subseteq \bigcup \{E_i : \chi_i(r) = 1\},
\end{equation}
which implies that
\begin{equation}
        H^{\alpha - 1}_{con}(A \cap f^{-1}(r))
         \le \sum \{(\diam E_i)^{\alpha - 1} : \chi_i(r) = 1\} = h(r).
\end{equation}
Moreover, if $\delta > 0$ and $\diam E_i < \delta$ for every $i$, then
\begin{equation}
        H^{\alpha - 1}_\delta(A \cap f^{-1}(r)) \le h(r)
\end{equation}
for every $r \in {\bf R}$.

\section{Level sets, 2}
\label{level sets, 2}
\setcounter{equation}{0}

        Let us continue with the same notations as in the previous
section, and suppose also now that
\begin{equation}
        H^\alpha(A) < \infty.
\end{equation}
For each $j \ge 1$, let $\{E_{i, j}\}_i$ be a covering of $A$ by
finitely or countably many subsets of $M$ such that
\begin{equation}
        \diam E_{i, j} < 1/j
\end{equation}
for every $i$, and
\begin{equation}
        \sum_i (\diam E_{i, j})^\alpha < H^\alpha_{1/j}(A) + \frac{1}{j}.
\end{equation}
Let $h_j(r)$ be the function on ${\bf R}$ corresponding to this covering
as in the previous section.  Thus
\begin{equation}
 \int_{\bf R} h_j(r) \, dr \le k \Big(H^\alpha_{1/j}(A) + \frac{1}{j}\Big)
\end{equation}
and
\begin{equation}
        H^{\alpha - 1}_{1/j}(A \cap f^{-1}(r)) \le h_j(r)
\end{equation}
for every $r \in {\bf R}$ and $j \ge 1$.  The latter implies that
\begin{equation}
        H^{\alpha - 1}(A \cap f^{-1}(r)) \le \liminf_{j \to \infty} h_j(r)
\end{equation}
for every $r \in {\bf R}$.  By Fatou's lemma,
\begin{equation}
 \int_{\bf R} \liminf_{j \to \infty} h_j(r) \, dr \le k \, H^\alpha(A).
\end{equation}
In particular, if $H^\alpha(A) = 0$, then it follows that
\begin{equation}
        H^{\alpha - 1}(A \cap f^{-1}(r)) = 0
\end{equation}
for almost every $r \in {\bf R}$.

\section{Topological dimension $\le n$}
\label{topological dimension le n}
\setcounter{equation}{0}

        A separable metric space $M$ is said to have \emph{topological
dimension $\le n$} if for every $p \in M$ and $r > 0$ there is an open
set $U$ in $M$ such that $p \in U \subseteq B(p, r)$ and $\partial U$
has topological dimension $\le n - 1$.  Thus the definition proceeds
inductively, and one can include the $n = 0$ case by considering the
empty set to have topological dimension $-1$.  More precisely, this
definition also applies to any subset $X$ of $M$, using the
restriction of the metric on $M$ to $X$.  If $M$ has topological
dimension $\le n$, then one can check that every $X \subseteq M$ has
topological dimension $\le n$, using induction.

        If $H^{n + 1}(M) = 0$ and $p \in M$, then the discussion in
the previous section applied to $\alpha = n + 1$ and $f_p(x) = d(p,
x)$ implies that
\begin{equation}
        H^n(\{x \in M : d(p, x) = r\}) = 0
\end{equation}
for almost every $r > 0$.  This implies that $M$ has topological
dimension $\le n$, because $\partial B(p, r)$ has topological
dimension $\le n - 1$ for every $p \in M$ and almost every $r > 0$, by
induction.  In particular, the topological dimension of $M$ is less
than or equal to the Hausdorff dimension of $M$, since $H^\alpha(M) =
0$ when $\alpha > {\dim}_H M$.

        A well-known theorem states that a separable metric space has
topological dimension $\le n$ when it is the union of a sequence of
closed subsets with topological dimension $\le n$, extending the $n =
0$ case discussed in Section \ref{closed subsets}.  As a consequence,
a separable metric space $M$ has topological dimension $\le n$ when $M
= X \cup Y$, $X$ and $Y$ have topological dimension $\le n$, and $X$
is closed.  For $M \backslash X$ is an open set with topological
dimension $\le n$ in this case, and one can argue as in Section
\ref{closed subsets} that $M$ is a countable union of closed sets with
topological dimension $\le n$.

        If $X$, $Y$ are arbitrary subsets of a separable metric space
with topological dimensions $\le m, n$, respectively, then it can be
shown that
\begin{equation}
\label{the topological dimension of X cup Y is le m + n + 1}
 \hbox{the topological dimension of $X \cup Y$ is $\le m + n + 1$.}
\end{equation}
This estimate is sharp, as in the example of the sets of rational and
irrational numbers.

        Another well-known theorem asserts that a separable metric
space $M$ with topological dimension $\le n$ can be expressed as the
union of two subsets $X$ and $Y$, where $X$ has topological dimension
$\le n - 1$, and $Y$ has topological dimension $0$.  Under these
conditions, there is a countable base for the topology of $M$
consisting of open sets whose boundaries have topological dimension
$\le n - 1$.  If $X$ is the union of the boundaries of the open sets
in this countable base, then $X$ has topological dimension $\le n -
1$, since it is the countable union of closed sets with topological
dimension $\le n - 1$.  It is easy to see that $Y = M \backslash X$
has topological dimension $0$, by construction.  Note that the theorem
about countable unions of closed sets is applied to sets of
topological dimension $\le n - 1$ in this argument, which permits the
conclusion to be used in the analysis of countable unions of closed
sets of topological dimension $\le n$, by induction.

        By repeating the process, $M$ can be realized as the union of
$n + 1$ subsets of topological dimension $0$.  This is consistent with
the estimate (\ref{the topological dimension of X cup Y is le m + n +
1}) for the topological dimension of the union of arbitrary subsets of
$M$.  The latter can also be used to show that $M$ is the union of a
pair of sets of topological dimension $\le l, r$, respectively, when
$n = l + r +1$.

        If $H^{n + 1}(M) = 0$, then there is a set $X \subseteq M$
such that $H^n(X) = 0$ and $M \backslash X$ has topological dimension
$0$.  This simply uses countable subadditivity of Hausdorff measure
instead of the theorem about topological dimensions of countable
unions of closed sets.

\section{Extrinsic conditions, 2}
\label{extrinsic conditions, 2}
\setcounter{equation}{0}

        Let $M$ be a metric space, let $X$ be a subset of $M$, and let
$E$ be a relatively open subset of $X$.  Of course, $X \backslash
\overline{E}$ is a relatively open set in $X$ that is disjoint from
$E$.  This implies that $E$ and $X \backslash \overline{E}$ are
separated as subsets of $X$ and hence in $M$, in the sense that
neither contains a limit point of the other.  By the strong version of
normality for metric spaces, there are disjoint open sets $U$, $V$ in
$M$ such that
\begin{equation}
        E \subseteq U, \ X \backslash \overline{E} \subseteq V.
\end{equation}
Let $\partial_X E$ be the boudary of $E$ relative to $X$, consisting
of limit points of $E$ in $X$ that are not contained in $E$.  It is
easy to see that
\begin{equation}
        \partial U \cap X \subseteq \partial_X E,
\end{equation}
since $U$ is contained in the closed set $M \backslash V$ and hence
$\overline{U} \subseteq M \backslash V$.  If $W$ is another open set
in $M$ such that $E \subseteq W \subseteq U$, then we also have that
\begin{equation}
        \partial W \cap X \subseteq \partial_X E.
\end{equation}

        Suppose now that $M$ is separable and $X$ has topological
dimension $\le n$.  Using the remarks in the preceding paragraph, one
can check that each $p \in X$ has arbitrarily small neighborhoods $W$
in $M$ such that $\partial W \cap X$ has topological dimension $\le n
- 1$.  This actually works for every $p \in M$, because $X \cup \{p\}$
also has topological dimension $\le n$.  For every $p \in M$ and $r > 0$,
\begin{equation}
        \partial B(p, r) \subseteq \{x \in M : d(p, x) = r\}.
\end{equation}
If $H^{n + 1}(X) = 0$, then $H^n(\partial B(p, r) \cap X) = 0$ for
almost every $r > 0$, and one can take $W = B(p, r)$ for such $r$ in
this case.

        Note that a separable metric space $M$ has topological
dimension $\le n$ if and only if for every $p \in M$ and closed set $B
\subseteq M$ with $p \not\in B$ there is an open set $U$ in $M$ such
that $p \in U$, $\overline{U} \cap B = \emptyset$, and $\partial U$
has topological dimension $\le n - 1$.  A well-known theorem states
that if $M$ has topological dimension $\le n$ and $A$, $B$ are
disjoint closed subsets of $M$, then there is an open set $U$ in $M$
such that $A \subseteq U$, $\overline{U} \cap B = \emptyset$, and
$\partial U$ has topological dimension $\le n - 1$.  The extrinsic
version of this theorem asserts that if $X \subseteq M$ has
topological dimension $\le n$ and $A$, $B$ are disjoint closed subsets
of $M$, then there is an open set $U \subseteq M$ such that $A
\subseteq U$, $\overline{U} \cap B = \emptyset$, and $\partial U \cap
X$ has topological dimension $\le n - 1$.  If $X \subseteq M$
satisfies $H^{n + 1}(X) = 0$ and $A$, $B$ are disjoint nonempty closed
subsets of $M$, then one can use the function $\phi$ in Section
\ref{disjoint closed sets} to show that there is an open set $U
\subseteq M$ such that $A \subseteq U$, $\overline{U} \cap B =
\emptyset$, and $H^n(\partial U \cap X) = 0$.

\section{Intersections}
\label{intersections}
\setcounter{equation}{0}

        Let $M$ be a separable metric space with topological dimension
$\le n - 1$, and let $A_1, B_1, \ldots, A_n, B_n$ be $n$ disjoint
pairs of closed subsets of $M$.  By repeating the extrinsic separation
theorem mentioned in the previous section, it is easy to see that
there are $n$ open subsets $U_1, \ldots, U_n$ of $M$ such that $A_i
\subseteq U_i$ and $\overline{U_i} \cap B_i = \emptyset$ for each $i$, and
\begin{equation}
 \partial U_1 \cap \partial U_2 \cap \cdots \cap \partial U_n = \emptyset.
\end{equation}
Using the Brouwer fixed-point theorem, one can then show that ${\bf
R}^n$ does not have topological dimension $\le n - 1$.  Of course,
${\bf R}^n$ does have topological dimension $\le n$.

        Here is a slightly different way to look at the case where
$H^n(M) = 0$.  Let $\phi_i$ be the function associated to the pair
$A_i$, $B_i$ as in Section \ref{disjoint closed sets} for $i = 1,
\ldots, n$, and let
\begin{equation}
        \Phi = (\phi_1, \ldots, \phi_n)
\end{equation}
be the combined mapping from $M$ into ${\bf R}^n$.  If
$\Lambda^i_\epsilon \subseteq M$ is associated to $A_i$, $B_i$ as
before, then the restriction of $\Phi$ to
\begin{equation}
        \Lambda^1_\epsilon \cap \cdots \cap \Lambda^n_\epsilon
\end{equation}
is Lipschitz for each $\epsilon > 0$.  Of course,
\begin{equation}
        H^n(\Lambda^1_\epsilon \cap \cdots \cap \Lambda^n_\epsilon)
         \le H^n(M) = 0
\end{equation}
for each $\epsilon > 0$, and hence
\begin{equation}
        H^n(\Phi(\Lambda^1_\epsilon \cap \cdots \cap \Lambda^n_\epsilon)) = 0.
\end{equation}
This implies that
\begin{equation}
        H^n(\Phi(M)) \le \sum_{j = 1}^\infty
         H^n(\Phi(\Lambda^1_{1/j} \cap \cdots \cap \Lambda^n_{1/j})) = 0.
\end{equation}
Since the open unit cube $(0, 1)^n$ has positive $n$-dimensional
Hausdorff measure, we get that
\begin{equation}
        (0, 1)^n \not\subseteq \Phi(M).
\end{equation}
If $r = (r_1, \ldots, r_n) \in (0, 1)^n \backslash \Phi(M)$, then
\begin{equation}
        U_i = \{x \in M : \phi_i(x) < r_i\}
\end{equation}
has the properties mentioned in the previous paragraph.

\section{Embeddings}
\label{embeddings}
\setcounter{equation}{0}

        Let $M$ be a separable metric space with topological dimension
$\le n$.  A famous theorem states that $M$ is homeomorphic to a
bounded set in ${\bf R}^{2n + 1}$.  Moreover, this set in ${\bf R}^{2n
+ 1}$ may be taken to have the property that its closure has Hausdorff
dimension $\le n$.  In particular, there is a metric on $M$ that
determines the same topology and with respect to which $M$ has
Hausdorff dimension $\le n$.  Of course, the proof of this theorem
relies heavily on the theory of the topological dimension, and so
cannot be used to derive basic results about the topological dimension
from properties of Hausdorff measure.  In some situations, changing
the metric on one part may not say much about the rest.  For example,
this theorem implies that there is a topologically-equivalent metric
on the set of irrational numbers with Hausdorff dimension $0$, and
$1$-dimensional Hausdorff measure $0$ in particular.  The set of
rational numbers already has these properties with respect to the
standard metric, but the real line still has topological dimension
$1$.

        If $M$ is already embedded in some ${\bf R}^l$, then it may
not be possible to deform this embedding using global homeomorphisms
on ${\bf R}^l$ to one with Hausdorff dimension equal to the
topological dimension as in the previous paragraph, even when $M$ is
compact.  See \cite{re, tr1, v1} for more information.

\section{Local Lipschitz conditions}
\label{locally lipschitz}
\setcounter{equation}{0}

        Let $(M_1, d_1(x, y))$ and $(M_2, d_2(u, v))$ be metric
spaces.  Let us say that a mapping $f : M_1 \to M_2$ is \emph{locally
$k$-Lipschitz at scale $\delta$} for some $k \ge 0$, $\delta > 0$ if
\begin{equation}
        d_2(f(x), f(y)) \le k \, d_1(x, y)
\end{equation}
for every $x, y \in M_1$ such that
\begin{equation}
        d_1(x, y) < \delta.
\end{equation}
As in Section \ref{lipschitz mappings}, if $f$ is locally $k$-Lipschitz
at the scale of $\delta$, $A \subseteq M_1$, and $\alpha \ge 0$, then
\begin{equation}
        H^\alpha_{k \delta'}(f(A)) \le k^\alpha \, H^\alpha_{\delta'}(A)
\end{equation}
when $0 < \delta' < \delta$, and hence
\begin{equation}
        H^\alpha(f(A)) \le k^\alpha \, H^\alpha(A).
\end{equation}
There are versions of the statements in Sections \ref{level sets} and
\ref{level sets, 2} for local Lipschitz conditions as well.

        For example, suppose that $M$ is a compact connected $C^1$
submanifold of ${\bf R}^n$.  Let $d_1(x, y)$ be the ordinary Euclidean
distance restricted to $x, y \in M$, and let $d_2(x, y)$ be the
Riemannian distance on $M$, which is to say the length of the shortest
path on $M$ connecting $x$ and $y$.  The latter is always greater than
or equal to the former, so that the identity mapping on $M$ is
$1$-Lipschitz as a mapping from $(M, d_2(x, y))$ onto $(M, d_1(x,
y))$.  The identity mapping on $M$ is also Lipschitz as a mapping from
$(M, d_1(x, y))$ onto $(M, d_2(x, y))$, and more precisely it is
locally Lipschitz at the scale of $\delta$ with a constant $k(\delta)
\to 1$ as $\delta \to 0$.  It follows that Hausdorff measures on $M$
with respect to these two metrics are the same.

        A mapping between metric spaces may be described as
\emph{locally Lipschitz} if each point in the domain has a
neighborhood on which the mapping is Lipschitz.  This permits both the
size of the neighborhood and the Lipschitz constant to depend on the
point. It is easy to see that the restriction of a locally Lipschitz
mapping to a compact set is Lipschitz.  Similarly, a mapping is
\emph{locally $k$-Lipschitz} for some $k \ge 0$ if each element of the
domain has a neighborhood on which the mapping is $k$-Lipschitz.  If
$f$ is locally $k$-Lipschitz and $A \subseteq M_1$ is compact, then
$f$ is locally $k$-Lipschitz at the scale of $\delta$ on $A$ for some
$\delta > 0$.  More precisely, one can cover $A$ by open balls $B(p,
r)$ such that the restriction of $f$ to $B(p, 2 \, r)$ is
$k$-Lipschitz.  Compactness implies that $A$ can be covered by
finitely many such balls $B(p_1, r_1), \ldots, B(p_n, r_n)$, and one
can take $\delta = \min(r_1, \ldots, r_n)$.

\section{Other Lipschitz conditions}
\label{other lipschitz conditions}
\setcounter{equation}{0}

        Let $(M_1, d_1(x, y))$, $(M_2, d_2(u, v))$ be metric spaces
again, and let $a$ be a positive real number.  A mapping $f : M_1 \to
M_2$ is said to be \emph{Lipschitz of order $a$} if there is a $k > 0$
such that
\begin{equation}
        d_2(f(x), f(y)) \le k \, d_1(x, y)^a
\end{equation}
for every $x, y \in M_1$.  Of course, this reduces to the ordinary
Lipschitz condition when $a = 1$.  This condition is also known as
\emph{H\"older continuity of order $a$}, and one can consider local
versions as well, as in the previous section.  In some situations,
such as for real-valued functions on Euclidean spaces, quite different
conditions are used for $a > 1$, related to the regularity of the
derivatives of $f$.

        If $a \le 1$, then $f$ is Lipschitz of order $a$ as a mapping
\begin{equation}
\label{(M_1, d_1(x, y)) to (M_2, d_2(u, v))}
        (M_1, d_1(x, y)) \to (M_2, d_2(u, v))
\end{equation}
if and only if it is Lipschitz of order $1$ as a mapping
\begin{equation}
        (M_1, d_1(x, y)^a) \to (M_2, d_2(u, v)).
\end{equation}
This works for any $a > 0$ when $d_1(x, y)$ is an ultrametric on
$M_1$, but otherwise $d_1(x, y)^a$ may not be a metric on $M_1$.  If
$d_1(x, y)$ is a snowflake metric, so that
\begin{equation}
        d_1(x, y) = d(x, y)^t
\end{equation}
for some metric $d(x, y)$ on $M$ and $t \in (0, 1)$, then this works
for $a \le 1/t$.  Alternatively, $f$ is Lipschitz of order $a \ge 1$
with respect to the initial metrics as in (\ref{(M_1, d_1(x, y)) to
(M_2, d_2(u, v))}) if and only if it is Lipschitz of order $1$ as a
mapping
\begin{equation}
        (M_1, d_1(x, y)) \to (M_2, d_2(u, v)^{1/a}).
\end{equation}
The same statement holds for every $a > 0$ when $d_2(u, v)$ is an
ultrametric, and for the appropriate range of $a$ when $d_2(u, v)$ is
a snowflake metric.

        In particular, the estimates for Hausdorff measures in Section
\ref{lipschitz mappings} carry over to Lipschitz mappings of any
order.  Specifically, if $f$ is $k$-Lipschitz of order $a$ and $E
\subseteq M_1$ is bounded, then $f(E)$ is bounded in $M_2$, and
\begin{equation}
        \diam f(E) \le k \, (\diam E)^a.
\end{equation}
This implies that
\begin{equation}
        H^{\alpha}(f(A)) \le k^\alpha \, H^{a \alpha}(A)
\end{equation}
for every $A \subseteq M_1$ and $\alpha \ge 0$, and hence
\begin{equation}
        {\dim}_H f(A) \le \frac{{\dim}_H A}{a}.
\end{equation}

\section{Quasimetrics}
\label{quasimetrics}
\setcounter{equation}{0}

        A \emph{quasimetric} on a set $M$ is a symmetric nonnegative
real-valued function $d(x, y)$ defined for $x, y \in M$ such that
$d(x, y) = 0$ if and only if $x = y$ and
\begin{equation}
        d(x, z) \le C \, (d(x, y) + d(y, z))
\end{equation}
for some $C \ge 1$ and all $x, y, z \in M$.  Thus $d(x, y)$ is a
metric on $M$ when $C = 1$.  If $d(x, y)$ is a quasimetric on $M$,
then
\begin{equation}
\label{d(x, y)^t}
        d(x, y)^t
\end{equation}
is also a quasimetric for every $t > 0$.  In particular, $d(x, y)^t$
is a quasimetric when $d(x, y)$ is a metric and $t > 1$.  Lipschitz
conditions can be extended to mappings between sets equipped with
quasimetrics instead of metrics in the obvious way.  The order of a
Lipschitz condition can be changed by changing the quasimetrics on the
domain or the range, is in the previous section.  This is a bit
simpler for quasimetrics, since arbitrary $t > 0$ are allowed in
(\ref{d(x, y)^t}).  However, an advantage of ordinary metrics is that
they determine real-valued Lipschitz functions of order $1$.  If $d(x,
y)$ is a quasimetric on $M$, then it is shown in \cite{m-s} that there
is a metric $\rho(x, y)$ on $M$ and $C_1, \delta > 0$ such that
\begin{equation}
        C_1^{-1} \, \rho(x, y) \le d(x, y)^\delta \le C_1 \, \rho(x, y)
\end{equation}
for every $x, y \in M$.

\section{Locally flat mappings}
\label{local flatness}
\setcounter{equation}{0}

        Let $(M_1, d_1(x, y))$, $(M_2, d_2(u, v))$ be metric spaces
again.  Let us say that a mapping $f : M_1 \to M_2$ is \emph{locally
flat} if it is locally Lipschitz, and for each $p \in M_1$ the
Lipschitz constant of the restriction of $f$ to $B(p, r)$ converges to
$0$ as $r \to 0$.  Equivalently, $f$ is locally flat if for each $p
\in M_1$ and $\epsilon > 0$ there is an $r > 0$ such that the
restriction of $f$ to $B(p, r)$ is $\epsilon$-Lipschitz, which is the
same as saying that $f$ is locally $\epsilon$-Lipschitz for each
$\epsilon > 0$.  Similarly, let us say that $f$ is \emph{uniformly
locally flat} if it is locally $k(\delta)$-Lipschitz at the scale of
$\delta$ for sufficiently small $\delta > 0$, where
\begin{equation}
\label{lim_{delta to 0} k(delta) = 0}
        \lim_{\delta \to 0} k(\delta) = 0.
\end{equation}
Thus uniformly locally flat mappings are locally flat, and the
restriction of a locally flat mapping to a compact set is uniformly
locally flat.

        A real-valued locally flat function on ${\bf R}^n$ has
differential equal to $0$ at each point.  If $f : M_1 \to M_2$ is
$k$-Lipschitz of order $a > 1$, then $f$ is uniformly locally flat,
with $k(\delta) = k \, \delta^{a - 1}$.  Of course, locally constant
mappings are locally flat, and may not be constant on disconnected
spaces.  Snowflake spaces can have nonconstant Lipschitz functions of
order $a > 1$ even when they are connected.

        If $f$ is uniformly locally flat, $A \subseteq M$, and
$H^\alpha(A) < + \infty$ for some $\alpha \ge 0$, then
\begin{equation}
        H^\alpha(f(A)) = 0.
\end{equation}
The same statement holds when $f$ is locally flat and $A$ is also
compact.  If $A$ is connected, then $f(A)$ is connected, and
\begin{equation}
        \diam f(A) \le H^1(f(A)),
\end{equation}
as in Section \ref{connected sets}.  Hence $f$ is constant on $A$ when
$A$ is connected and $H^1(A)$ is finite.  The same argument works for
Lipschitz mappings of order $a > 1$ when $A$ is connected and $H^a(A)
= 0$.

\section{Rectifiable curves}
\label{rectifiable curves}
\setcounter{equation}{0}

        Let $a$, $b$ be real numbers with $a \le b$, let $(M, d(u,
v))$ be a metric space, and let $p(t)$ be a continuous path in $M$
defined on the closed interval $[a, b]$.  If $\mathcal{P} = \{t_j\}_{j
= 0}^n$ is a partition of $[a, b]$, so that
\begin{equation}
        a = t_0 < t_1 < \cdots < t_n = b,
\end{equation}
then put
\begin{equation}
        \Lambda(\mathcal{P}) = \sum_{j = 1}^n d(p(t_j), p(t_{j - 1})).
\end{equation}
The path $p(t)$, $a \le t \le b$, is said to have \emph{finite length}
if the $\Lambda(\mathcal{P})$'s are bounded, in which case the length
of the path is defined by
\begin{equation}
        \Lambda_a^b = \sup_{\mathcal{P}} \Lambda(\mathcal{P}),
\end{equation}
where the supremum is taken over all partitions $\mathcal{P}$ of $[a, b]$.

        For example, if $p : [a, b] \to M$ is $k$-Lipschitz, then
\begin{equation}
        \Lambda(\mathcal{P}) \le k \, (b - a)
\end{equation}
for every partition $\mathcal{P}$ of $[a, b]$, and hence $p$ has
finite length
\begin{equation}
        \Lambda_a^b \le k \, (b - a).
\end{equation}
In particular, constant paths have length $0$.  Conversely, for any path $p$,
\begin{equation}
        d(p(x), p(y)) \le \Lambda(\mathcal{P})
\end{equation}
when $a \le x \le y \le b$ and $\mathcal{P}$ is the partition
consisting of these four points, which implies that
\begin{equation}
        \diam p([a, b]) \le \Lambda_a^b.
\end{equation}
This shows that paths with length $0$ are constant.

        If $a \le x \le y \le b$, then every partition of $[x, y]$ can
be extended to a partition of $[a, b]$.  Hence the restriction of a
continuous path $p : [a, b] \to M$ of finite length to $[x, y]$ also
has finite length, and
\begin{equation}
\label{Lambda_x^y le Lambda_a^b}
        \Lambda_x^y \le \Lambda_a^b.
\end{equation}
A partition $\mathcal{P} = \{t_j\}_{j = 0}^n$ of $[a, b]$ is said to
be a \emph{refinement} of another partition $\mathcal{P}' = \{r_i\}_{i
= 0}^l$ if for each $i = 0, 1, \ldots, l$ there is a $j = 0, 1,
\ldots, n$ such that $r_i = t_j$.  Using the triangle inequality, one
can check that
\begin{equation}
\label{Lambda(mathcal{P}') le Lambda(mathcal{P})}
        \Lambda(\mathcal{P}') \le \Lambda(\mathcal{P})
\end{equation}
when $\mathcal{P}$ is a refinement of $\mathcal{P}'$.  If $p$ has
finite length and $a \le x \le b$, then it follows that
\begin{equation}
        \Lambda_a^b = \Lambda_a^x + \Lambda_x^b.
\end{equation}
Indeed, $\Lambda_a^x + \Lambda_x^b \le \Lambda_a^b$ because arbitrary
partitions of $[a, x]$ and $[x, b]$ can be combined to get a partition
of $[a, b]$.  The opposite inequality holds because any partition of
$[a, b]$ can be refined to include $x$, and the refinement is then a
combination of partitions of $[a, x]$ and $[x, b]$.  Similar reasoning
shows that $p : [a, b] \to M$ has finite length when the restrictions
of $p$ to $[a, x]$ and $[x, b]$ have finite length.

        In totally disconnected spaces, continuous paths are
automatically constant.  Snowflake spaces may contain many nontrivial
continuous paths, but one can check that continuous paths of finite
length are constant.  However, there are fractal sets such as
Sierpinski gaskets and carpets and Menger sponges with numerous
nontrivial continuous paths of finite length.  Smooth manifolds also
have plenty of continuous paths of finite length, and this can be
extended to sub-Riemannian spaces as well.

\section{Length and measure}
\label{length and measure}
\setcounter{equation}{0}

        Let $(M, d(u, v))$ be a metric space, and let $p : [a, b] \to
M$ be a continuous path of finite length. If $\mathcal{P} = \{t_j\}_{j
= 0}^n$ is any partition of $[a, b]$, then
\begin{equation}
        \sum_{j = 1}^n \diam p([t_{j - 1}, t_j])
         \le \sum_{j = 1}^n \Lambda_{t_{j - 1}}^{t_j} = \Lambda_a^b,
\end{equation}
by the computations in the previous section.  If in addition $\delta > 0$ and
\begin{equation}
\label{diam p([t_{j - 1}, t_j]) < delta, 1 le j le n}
        \diam p([t_{j - 1}, t_j]) < \delta, \ 1 \le j \le n,
\end{equation}
then it follows that
\begin{equation}
        H^1_\delta(p([a, b])) \le \Lambda_a^b.
\end{equation}
Because of uniform continuity, (\ref{diam p([t_{j - 1}, t_j]) < delta,
1 le j le n}) holds for any $\delta > 0$ when $\mathcal{P}$ is a
sufficiently fine partition of $[a, b]$, which implies that
\begin{equation}
\label{H^1(p([a, b])) le Lambda_a^b}
        H^1(p([a, b])) \le \Lambda_a^b.
\end{equation}

        If $p : [a, b] \to M$ is also injective, then
\begin{equation}
        H^1(p([a, b])) = \Lambda_a^b.
\end{equation}
To see this, observe that
\begin{equation}
        d(p(x), p(y)) \le \diam p([x, y]) \le H^1(p([x, y]))
\end{equation}
when $a \le x \le y \le b$, since $p([x, y])$ is connected.  If $I_1,
\ldots, I_l$ are pairwise disjoint closed subintervals of $[a, b]$,
then $p(I_1), \ldots, p(I_l)$ are pairwise disjoint compact subsets of
$M$, and are therefore at positive distance from each other.  Thus
\begin{equation}
        \sum_{i = 1}^l H^1(p(I_i)) = H^1\Big(\bigcup_{i = 1}^l p(I_i)\Big)
                                    \le H^1(p([a, b])),
\end{equation}
using the additivity of Hausdorff measure in this case, as in Section
\ref{hausdorff measures}.  If
\begin{equation}
        a \le x_1 < y_1 < x_2 < y_2 < x_3 < \cdots < y_l \le b,
\end{equation}
then we can apply these estimates to $I_i = [x_i, y_i]$ to get that
\begin{equation}
        \sum_{i = 1}^l d(p(x_i), p(y_i)) \le H^1(p([a, b])).
\end{equation}
This also holds when
\begin{equation}
        a \le x_1 < y_1 \le x_2 < y_2 \le x_3 < \cdots < y_l \le b,
\end{equation}
by passing to suitable limits.  Hence
\begin{equation}
        \Lambda(\mathcal{P}) \le H^1(p([a, b]))
\end{equation}
for every partition $\mathcal{P}$ of $[a, b]$, as desired.  Note that
an injective continuous path $p : [a, b] \to M$ has finite length when
$H^1(p([a, b])) < +\infty$, by the same argument.

        This argument can also be extended to deal with paths with
only finitely many crossings, for instance.  However, strict
inequality can occur in (\ref{H^1(p([a, b])) le Lambda_a^b}) for
arbitrary continuous paths of finite length, as when such a path
retraces an arc.

\section{Mappings of paths}
\label{mappings of paths}
\setcounter{equation}{0}

        Let $(M_1, d_1(u, v))$ and $(M_2, d_2(w, z))$ be metric
spaces, and let $p : [a, b] \to M_1$ and $f : M_1 \to M_2$ be
continuous mappings.  Thus $\widetilde{p} = f \circ p : [a, b] \to
M_2$ is also continuous.  For each partition $\mathcal{P}$ of $[a,
b]$, let $\Lambda(\mathcal{P})$, $\widetilde{\Lambda}(\mathcal{P})$ be
the corresponding sums for the paths $p$, $\widetilde{p}$,
respectively, as in Section \ref{rectifiable curves}.

        If $f : M_1 \to M_2$ is $k$-Lipschitz for some $k \ge 0$, then
\begin{equation}
        \widetilde{\Lambda}(\mathcal{P}) \le k \, \Lambda(\mathcal{P})
\end{equation}
for every partition $\mathcal{P}$ of $[a, b]$.  Hence $\widetilde{p}$
has finite length when $p$ does, and
\begin{equation}
        \widetilde{\Lambda}_a^b \le k \, \Lambda_a^b,
\end{equation}
where $\Lambda_a^b$, $\widetilde{\Lambda}_a^b$ are the lengths of $p$,
$\widetilde{p}$, respectively.

        If $f : M_1 \to M_2$ is locally $k$-Lipschitz at the scale of
$\delta$ for some $\delta > 0$, then we get the same estimate for
$\widetilde{\Lambda}(\mathcal{P})$ when the partition $\mathcal{P} =
\{t_j\}_{j = 0}^n$ is sufficiently fine so that
\begin{equation}
        d_1(p(t_j), p(t_{j - 1})) < \delta, \ j = 1, \ldots, n.
\end{equation}
Every partition of $[a, b]$ has a refinement with this property,
because of uniform continuity.  It follows again that $\widetilde{p}$
has finite length when $p$ does, and with the same estimate for the
length.

        If $f : M_1 \to M_2$ is locally $k$-Lipschitz, then the
restriction of $f$ to $p([a, b])$ is locally $k$-Lipschitz at the
scale of $\delta$ for some $\delta > 0$, since $p([a, b])$ is compact.
Hence the same conclusions hold in this case.  If $f : M_1 \to M_2$ is
locally flat, then $\widetilde{p}$ has length $0$, and is therefore constant.

\section{Reparameterizations}
\label{reparameterizations}
\setcounter{equation}{0}

        Let $(M, d(u, v))$ be a metric space, and let $p : [a, b] \to
M$ be a continuous mapping.  Also let $\phi$ be a strictly increasing
continuous mapping from another closed interval $[\widehat{a},
\widehat{b}]$ onto $[a, b]$, so that $\widehat{p} = p \circ \phi :
[\widehat{a}, \widehat{b}] \to M$ is a continuous path in $M$.  
If $\widehat{P} = \{\widehat{t}_j\}_{j = 0}^n$ is a partition of
$[\widehat{a}, \widehat{b}]$, then $\mathcal{P} =
\{\phi(\widehat{t}_j)\}_{j = 0}^n$ is a partition of $[a, b]$, and
every partition of $[a, b]$ corresponds to a partition $[\widehat{a},
\widehat{b}]$ in this way.  By construction,
\begin{equation}
 \sum_{j = 1}^n d(\widehat{p}(\widehat{t}_j), \widehat{p}(\widehat{t}_{j - 1}))
  = \sum_{j = 1}^n d(p(\phi(\widehat{t}_j)), p(\phi(\widehat{t}_{j - 1}))),
\end{equation}
which implies that $\widehat{p}$ has finite length if and only if $p$
has finite length, and that their lengths are the same in this case.

        The same conclusion holds when $\phi$ is a monotone increasing
continuous mapping from $[\widehat{a}, \widehat{b}]$ onto $[a, b]$.
The correspondence between partitions of $[\widehat{a}, \widehat{b}]$
and $[a, b]$ is a bit more complicated when $\phi$ is not strictly
increasing, but this is not significant for the approximations of the
lengths of $\widehat{p}$ and $p$.  This is because $\phi$ and hence
$\widehat{p}$ is constant on any interval $[x, y] \subseteq
[\widehat{a}, \widehat{b}]$ such that $\phi(x) = \phi(y)$.

        If $p : [a, b] \to M$ has finite length, then the length
$\Lambda_a^r$ of the restriction of $p$ to $[a, r]$ is a monotone
increasing function of $r$ on $[a, b]$.  One can also show that
$\Lambda_a^r$ is continuous in $r$, as follows.  Let $\epsilon > 0$ be
given, and let $\mathcal{P}$ be a partition of $[a, b]$ such that
\begin{equation}
        \Lambda_a^b < \Lambda(\mathcal{P}) + \epsilon.
\end{equation}
Suppose that $a \le r < t \le b$ and that $\mathcal{P}$ does not
contain any element of $(r, t)$.  Let $\mathcal{P}(r, t)$ be a
refinement of $\mathcal{P}$ that contains $r$, $t$ as consecutive
terms.  Thus
\begin{equation}
        \Lambda(P) \le \Lambda(\mathcal{P}(r, t)) 
                    \le \Lambda_a^r + d(p(r), p(t)) + \Lambda_t^b,
\end{equation}
since the rest of $\mathcal{P}(r, t)$ partitions $[a, r]$ and $[t,
b]$.  This implies that
\begin{equation}
        \Lambda_r^t \le d(p(r), p(t)) + \epsilon,
\end{equation}
because $\Lambda_a^b = \Lambda_a^r + \Lambda_r^t + \Lambda_t^b$.  This
permits $\Lambda_r^t$ to be estimated in terms of the continuity of
$p$ when $r$, $t$ are sufficiently close.

        If $a \le r < t \le b$ and $\Lambda_a^r = \Lambda_a^t$, then
$\Lambda_r^t = 0$, and $p$ is constant on $[r, t]$.  It follows that
there is a mapping $q : [0, \Lambda_a^b] \to M$ such that $p(r) =
q(\Lambda_a^r)$ for each $r \in [a, b]$.  Moreover, $q$ is
$1$-Lipschitz, because $d(p(r), p(t)) \le \Lambda_r^t$ when $a \le r <
t \le b$.  The length of $q$ is equal to the length $\Lambda_a^b$ of
$p$, by the earlier remarks about arbitrary reparameterizations.

\end{document}